\documentclass[a4paper,12pt]{article}
\usepackage{amssymb,latexsym,amsmath,amsthm,amscd,graphpap,color}
\usepackage{amsmath}
\usepackage{amsfonts}
\usepackage{amssymb}
\usepackage[dvips]{epsfig}
\usepackage{verbatim}
\newtheorem{theorem}{Theorem}[section]
\newtheorem{prop}[theorem]{Proposition}

\newtheorem{lemma}[theorem]{Lemma}
\newtheorem{corollary}[theorem]{Corollary}
\newenvironment{demo}{ \noindent \emph{\textbf{Proof:}}}{\hfill$\square$\\
\vspace{0.4cm}}

\newcommand{\RR}{\mathbb{R}}

\newcommand{\Lc}{\mathcal{L}}

\newcommand{\NN}{\mathbb{N}}

\newcommand{\Gg}{\mathfrak{G}}
\newcommand{\Cg}{\mathfrak{C}}

\newcommand{\Cc}{\mathcal{C}}
\newcommand{\Nc}{\mathcal{N}}

\newcommand{\Oc}{\mathcal{O}}

\newcommand{\Ac}{\mathcal{A}}

\newcommand{\Kc}{\mathcal{K}}

\newcommand{\Uc}{\mathcal{U}}

\newcommand{\Drond}[2]{\frac{\partial #1}{\partial #2}}

\newcommand{\supp}{\mathrm{supp}}
\newcommand{\no}{n$^{\text{o}}$}
\newcommand{\grad}{\nabla}
\renewcommand{\div}{\mathrm{div}}
\renewcommand{\th}{\mathrm{th}}

\newcommand{\pc}{ \usefont{T1}{cmtl}{m}{n} \selectfont}
\newcounter{compte}
\newenvironment{enum2}{\begin{list}{\roman{compte})} {\usecounter{compte}
\topsep=1mm \itemsep=0.2mm \leftmargin=5mm  } }
{\end{list}}

\numberwithin{equation}{section}

\newdimen\texpscorrection
\newdimen\figcenter
\def\figurewithtex #1 #2 #3 #4 #5\cr{\null
  {\goodbreak\figcenter=\hsize\relax
  \advance\figcenter by -#4truecm
  \divide\figcenter by 2
  \begin{figure}[hbt]
  \vskip #3truecm\noindent\hskip\figcenter
  \includegraphics{#1}{\hskip\texpscorrection\input #2 }
  \vskip 0.8truecm{\baselineskip=0.8\baselineskip
  \noindent \vbox{\noindent {\footnotesize #5}}\par}
  \end{figure}}}
\def\point#1 #2 #3 {\rlap{\kern #1 truecm
\raise #2 truecm \hbox{#3}}}

\usepackage[english,activeacute]{babel}

\usepackage[latin1]{inputenc}
\usepackage{verbatim}

\newcommand{\R}{\mathbb{R}}

\newcommand\nor[2]{\left\|#1\right\|_{#2}}%la norme
\newcommand\bna{\begin{eqnarray*}}%�quation non num�rot�e
\newcommand\ena{\end{eqnarray*}}
\newcommand\bnan{\begin{eqnarray}}%�quation num�rot�e
\newcommand\enan{\end{eqnarray}}

\newcommand\bneq{\begin{eqnarray*}\left\lbrace \begin{array}{rcl}}%systeme equation sans numerotation
\newcommand\eneq{\end{array} \right.\end{eqnarray*}}
\newcommand\bneqn{\begin{eqnarray}\left\lbrace \begin{array}{rcl}}%systeme equation avec numerotation
\newcommand\eneqn{\end{array} \right.\end{eqnarray}}

\newcommand\bnp{\begin{demo}}
\newcommand\enp{\end{demo}}

\newcommand\e{\varepsilon}

\begin{document}

\title{\bf Stabilization for the semilinear wave equation with geometric
control condition}

\author{ Romain \textsc{Joly}\footnote{Institut Fourier - UMR5582
CNRS/Universit\'e de Grenoble - 100, rue des Maths - BP 74 - F-38402
St-Martin-d'H\`eres, France, email: {\pc romain.joly@ujf-grenoble.fr}} {~\&~}
Camille \textsc{Laurent}\footnote{CNRS, UMR 7598, Laboratoire Jacques-Louis
Lions, F-75005, Paris, France  }
\footnote{UPMC Univ Paris 06, UMR 7598, Laboratoire Jacques-Louis Lions,
F-75005, Paris, France, email: {\pc laurent@ann.jussieu.fr } }}
\maketitle
%\tableofcontents
%\vspace{1cm}

\begin{abstract}
In this article, we prove the exponential stabilization of the semilinear wave
equation with a damping effective in a zone satisfying the geometric
control condition only. The nonlinearity is assumed to be subcritical, defocusing and analytic.
The main novelty compared to previous results, is the proof of a unique
continuation result in large time for some undamped equation. The idea is to use
an asymptotic smoothing effect proved by Hale and Raugel in the context of
dynamical systems. Then, once the analyticity in time is proved, we apply a
unique continuation result with partial analyticity due to
Robbiano, Zuily, Tataru and H\"ormander. Some other consequences are also given
for the controllability and the existence of a compact attractor.\\[3mm]
{\sc Key words:} damped wave equation, stabilization,
analyticity, unique continuation property, compact attractor.\\
{\sc AMS subject classification:} 35B40, 35B60, 35B65,
35L71, 93D20, 35B41.\\
\end{abstract}
\renewcommand{\abstractname}{R\'esum\'e}
\begin{abstract}
Dans cet article, on prouve la d\'ecroissance exponentielle de l'\'equation des
ondes semilin\'eaires avec un amortissement actif dans
une zone satisfaisant seulement la condition de contr\^ole g\'eom\'etrique.
La
nonlin\'earit\'e est suppos\'ee sous-critique, d\'efocalisante et analytique. La
principale
nouveaut\'e par rapport aux r\'e\-sul\-tats pr\'e\-c\'e\-dents est la preuve
d'un
r\'esultat de prolongement unique en grand temps pour une solution non amortie.
L'id\'ee est d'utiliser un effet r\'egularisant asymptotique prouv\'e par Hale
et Raugel dans le contexte des syst\`emes dynamiques. Ensuite, une fois
l'analyticit\'e en temps prouv\'ee, on applique un th\'eor\`eme de prolongement
unique avec analyticit\'e partielle d\^u \`a Robbiano, Zuily, Tataru et
H\"ormander. Des applications \`a la contr\^olabilit\'e et \`a l'existence
d'attracteur global compact pour l'\'equation des ondes sont aussi
donn\'ees.\\[3mm]
\end{abstract}
\renewcommand{\abstractname}{Abstract}

\section{Introduction}
In this article, we consider
the semilinear damped wave equation
\begin{equation}\label{eq}
\left\{\begin{array}{ll}
\Box u +\gamma(x) \partial_t u+\beta u+f(u)=0~~&(t,x)\in\RR_+\times\Omega~,\\
u(t,x)=0&(t,x)\in\RR_+\times\partial\Omega\\
(u,\partial_t u)=(u_0,u_1)\in H^1_0(\Omega)\times L^2(\Omega)&
\end{array}\right.
\end{equation}
where $\Box=\partial_{tt}^2-\Delta$ with $\Delta$ being the Laplace-Beltrami
operator with Dirichlet boundary conditions. The domain $\Omega$ is a connected $\Cc^\infty$ three-dimensional
Riemannian manifold with boundaries, which is either:
\begin{enum2}
\item compact.
\item a compact perturbation of $\RR^3$, that is $\R^3\setminus D$ where $D$ is
a bounded smooth domain, endowed with a smooth metric equal to the euclidean one
outside of a ball.   
\item or a manifold with periodic geometry (cylinder, $\RR^3$ with periodic
metric etc.). 
\end{enum2}

The nonlinearity $f\in\Cc^1(\RR,\RR)$ is assumed to be defocusing, energy
subcritical and such that $0$ is an equilibrium point. More precisely, we assume
that there exists $C>0$  such that 
\begin{equation}\label{hyp-f}
f(0)=0~,~~sf(s)\geq 0~,~~|f(s)|\leq C(1+|s|)^p~\text{ and }~~|f'(s)|\leq
C(1+|s|)^{p-1}
\end{equation}
with $1\leq p < 5$. 

We assume $\beta\geq 0$ to be such that $\Delta-\beta$ is a negative definite
operator, that is that we have a Poincar\'e inequality $\int_\Omega |\grad
u|^2+\beta |u|^2 \geq C\int_\Omega |u|^2$ with $C>0$. In particular, it
may require $\beta>0$ if $\partial\Omega=\emptyset$ or if $\Omega$ is
unbounded. 
 
The damping $\gamma\in L^{\infty}(\Omega)$ is a non-negative function.
We assume that there exist an open set $\omega\subset\Omega$, $\alpha\in\RR$, $x_0\in\Omega$ and
$R\geq 0$ such that 
\begin{equation}\label{hyp-gamma}
\forall x \in\omega~,~\gamma(x)\geq\alpha>0~ \text{ and } \Omega\setminus B(x_0,R)
\subset \omega~.
\end{equation}
Moreover, we assume that $\omega$ satisfies the geometric control condition
introduced in \cite{RT} and \cite{BLR}
\begin{itemize}
\item[\bf (GCC)] There exists $L>0$ such that any generalized geodesic of
$\Omega$ of length $L$ meets the set $\omega$ where the damping is effective.
\end{itemize}
The associated energy
$E\in\Cc^0(X,\RR_+)$ is given by
\begin{equation}\label{def-E}
E(u):=E(u,\partial_t u)=\frac 12 \int_\Omega (|\partial_t u|^2+|\nabla
u|^2+\beta|u|^2) + \int_\Omega V(u)~,
\end{equation}
where $V(u)=\int_0^u f(s)ds$. Due to Assumption (\ref{hyp-f}) and the Sobolev
embedding $H^1(\Omega)\hookrightarrow L^6(\Omega)$, this energy is well defined
and moreover, if $u$ solves \eqref{eq}, we have, at least formally, 
\begin{equation}
\label{energydecay}
\partial_t E(u(t))=-\int_\Omega \gamma(x)|\partial_t u(x,t)|^2~dx~\leq~ 0~.
\end{equation}
The system is therefore dissipative. We are interested in the exponential decay
of the energy of the nonlinear damped wave equation \eqref{eq}, that is the
following property:
\begin{itemize}
 \item[\bf (ED)]For any $E_0\geq 0$, there exist $K>0$ and $\lambda>0$ such
that, for all solutions $u$ of \eqref{eq} with $E(u(0))\leq E_0$, 
\begin{equation*}
\forall t\geq 0~,~~E(u(t))\leq Ke^{-\lambda t} E(u(0))
\end{equation*}
\end{itemize}
Property (ED) means that the damping term $\gamma \partial_t u$ stabilizes any
solution of \eqref{eq} to zero, which is an important property from the
dynamical and control points of view.

\bigskip

\noindent Our main theorem is as follows.
\begin{theorem}
\label{th1}
Assume that the damping $\gamma$ satisfies \eqref{hyp-gamma} and the geometric
control condition (GCC). If $f$ is real analytic and satisfies \eqref{hyp-f},
then the exponential decay property (ED) holds.
\end{theorem}

Theorem \ref{th1} applies for nonlinearities $f$ which are globally analytic. 
Of course, the nonlinearities $f(u)=|u|^{p-1}u$ are not analytic if $p\not\in\{
1,3\}$, but we can replace these usual nonlinearities by similar ones
as $f(u)=(u/\th(u))^{p-1} u$, which are analytic for all $p\in[1,5)$. Note that
the estimates \eqref{hyp-f} are only required for $s\in \R$, so that it does not
imply that $f$ is polynomial.  
Moreover, we enhance that (ED) holds in fact for almost all the
nonlinearities $f$ satisfying \eqref{hyp-f}, including non-analytic ones.

More precisely, we set 
\begin{equation}\label{eq-Cg}
\Cg^1(\RR)=\{f\in\Cc^1(\RR) \text{ such that there exist }C>0\text{ and
}p\in[1,5)\text{ such that \eqref{hyp-f} holds }\}
\end{equation}
endowed with Whitney topology (or any other reasonable topology). We recall that
Whitney topology is the topology generated by the
neighbourhoods
\begin{equation}\label{eq-Whitney}
\Nc_{f,\delta}=\{~g\in\Cg^1(\RR)~|~\forall u\in\RR~,~~\max
(|f(u)-g(u)|,|f'(u)-g'(u)|) < \delta(u)~\}
\end{equation}
where $f$ is any function in $\Cg^1(\RR)$ and $\delta$ is any positive
continuous function. The set $\Cg^1(\RR)$ is a 
Baire space, which means that any generic set, that is any set
containing a countable intersection of open and dense sets, is dense in
$\Cg^1(\RR)$ (see Proposition \ref{prop-Baire}). Baire property ensures that
the genericity of a set in $\Cg^1(\RR)$ is a good notion for ``the set contains
almost all non-linearity $f$''. 
\begin{theorem}\label{th2}
Assume that the damping $\gamma$ satisfies \eqref{hyp-gamma} and the geometric
control condition (GCC). There exists a generic set $\Gg\subset \Cg^1(\RR)$ such
that the exponential decay property (ED) holds for all $f\in\Gg$.
\end{theorem}

\noindent The statements of both theorems lead to some remarks.\\
$\bullet$ Of course, our results and their proofs should easily extend to
any space dimension $d\geq 3$ if the exponent $p$ of the nonlinearity satisfies
$p<(d+2)/(d-2)$.\\
$\bullet$ Actually, it may be possible to get $\lambda>0$ in (ED) uniform with
respect to the size of the data. We can take for instance 
$\lambda=\tilde{\lambda}-\e$ where $\tilde{\lambda}$ is the decay rate of the
linear equation. The idea is that once we know the existence of a decay rate, we
know that the solution is close to zero for a large time. Then, for small
solutions, the nonlinear term can be neglected to get almost the same decay rate
as the linear equation. We refer for instance to \cite{KdV_LRZ} in the context
of KdV equation. Notice that the possibility to get the same result with a
constant $K$ independent on $E_0$ is an open problem. \\ 
$\bullet$ The assumption on $\beta$ is important to ensure some
coercivity of the energy and to preclude the spatially constant functions to be
undamped solutions for the linear equation. It has been proved in
\cite{DehPGNLW} for $\RR^3$ and in \cite{LaurentNLW} for a compact manifold that
exponential decay can fail without this term $\beta$.\\
$\bullet$ The geometric control condition is known to be not only
sufficient but also necessary for the exponential decay of the linear damped
equation. The proof of the optimality uses some sequences of solutions which are
asymptotically concentrated outside of the damping region. We can use the same
idea in our nonlinear stabilization context. First, the observability for a
certain time eventually large is known to be equivalent to the exponential decay
of the energy.  This was for instance noticed in \cite{DehPGNLW} Proposition 2,
in a similar context, see also Proposition \ref{critere-decay} of this paper.
Then, we take as initial data the same sequence that would 
give a counterexample for the linear observability. The linearizability property
(see \cite{linearisationondePG}) allows to obtain that the nonlinear solution is
asymptotically close to the linear one. This contradicts the observability
property for the nonlinear solution as it does for the linear case. Hence, the
geometric control condition is also necessary for the exponential decay of the
nonlinear equation.\\
$\bullet$ Our geometrical hypotheses on $\Omega$ may look strange, however
they are only assumed for sake of simplicity. In fact, our results should apply
more generally for any smooth manifold with bounded
geometry, that is that $\Omega$ can be covered by a set
of $\Cc^\infty-$charts $\alpha_i:U_i\longmapsto \alpha_i(U_i)\subset\RR^3$ such
that $\alpha_i(U_i)$ is equal either to $B(0,1)$ or to $B_+(0,1)=\{x\in B(0,1),
x_1>0\}$ (in the case with boundaries) and such that, for any
$r\geq 0$ and $s\in[1,\infty]$, the $W^{r,s}-$norm of a function $u$ in
$W^{r,s}(\Omega,\RR)$ is equivalent to the norm $(\sum_{i\in\NN}
\|u\circ\alpha_i^{-1}\|^s_{W^{r,s}(\alpha_i(U_i))})^{1/s}$.

\bigskip

The stabilization property (ED) for Equation \eqref{eq} has been studied in
\cite{Haraux2}, \cite{Zuazua1}, \cite{Zuazua2} and \cite{Dehman} for $p<3$. For
$p\in[3,5)$, our main reference is the work of Dehman, Lebeau and
Zuazua \cite{DLZ}. This work is mainly concerned with the
stabilization problem previously described, on the Euclidean space $\RR^3$ with
flat metric and stabilization active outside of a ball. The main purpose of this
paper is to extend their result to a non flat geometry where multiplier methods
cannot be used or do not give the optimal result with respect to the geometry.
Other stabilization results for the nonlinear wave equation can be found 
in \cite{AIN} and the references therein. Some works have been done in the
difficult critical case $p=5$, we refer to \cite{DehPGNLW} and
\cite{LaurentNLW}.

The proofs of these articles use three main ingredients:
\begin{enum2}
\item[(i)] the exponential decay of the linear equation, which is equivalent to
the geometric control condition (GCC),
\item[(ii)] a more or less involved compactness argument,
\item[(iii)] a unique continuation result implying that $u\equiv 0$ is the
unique solution of
\bneqn\label{eq-iii}
\Box u+\beta u+f(u)&=&0\\
\partial_tu&=&0\quad \textnormal{ on }[-T,T]\times \omega~. 
\eneqn
\end{enum2}
The results are mainly of the type: ``geometric control condition'' + ``unique
continuation'' implies ``exponential decay''. This type of implication is even
stated explicitly in some related works for the nonlinear Schr\"odinger equation
\cite{control-nl} and \cite{LaurentNLSdim3}.

In the subcritical case $p<5$, the less understood point is the unique
continuation property (iii). In the previous works as \cite{DLZ}, the authors
use unique continuation results based on Carleman estimates. The
resulting geometric assumptions are not very natural and are stronger than
(GCC). Indeed, the unique continuation was often proved with some Carleman
estimates that required some strong geometric conditions. For instance for a
flat metric, the usual geometric assumption that appear are often of
``multiplier type'' that is $\omega$ is a neighbourhood of $\left\{x\in \partial
\Omega \left|(x-x_0)\cdot n(x)>0\right.\right\}$ which are known to be stronger
than the geometric control condition (see \cite{Millerescape} for a discussion
about the links between these assumptions). Moreover, on curved spaces,
this type of condition often needs to be checked by hand in each situation,
which is mostly impossible.

Our main improvement in this paper is the proof of unique continuation in
infinite time under the geometric control condition only. We show  that,
if the nonlinearity $f$ is analytic (or generic), then one can use the result of
Robbiano and Zuily \cite{Rob-Zui} to obtain a unique continuation property (iii)
for infinite time $T=+\infty$ with the geometric control condition (GCC) only.

The central argument of the proof of our main result, Theorem \ref{th1}, is the
unique continuation property of \cite{Rob-Zui} (see Section
\ref{sect-UCP}). This result applies for solutions $u$ of \eqref{eq-iii}
being smooth in space and analytic in time. If $f$ is analytic, then the
solutions of \eqref{eq} are of course not necessarily
analytic in time since the damped wave equations are not smoothing in finite
time. However, the damped wave equations admit an asymptotic smoothing effect,
i.e. are smoothing in infinite time. Hale and Raugel have shown in
\cite{Hale-Raugel} that, for compact trajectories, this asymptotic smoothing effect also concerns the
analyticity (see Section \ref{sect-HR}). In other words, combining
\cite{Rob-Zui} and \cite{Hale-Raugel} shows that the unique solution of
\eqref{eq-iii} is $u\equiv 0$ if $f$ is analytic and if $T=+\infty$.
This combination has already been used by dynamicians in \cite{HR2} and
\cite{RJ} for $p<3$.

\bigskip

One of the main interests of this paper is the use of arguments coming from
both the dynamical study and the control theory of the damped wave equations.
The reader familiar with the control theory could find interesting the use of
the asymptotic smoothing effect to get unique continuation property with smooth
solutions. The one familiar with the dynamical study of PDEs could be interested
in the use of Strichartz estimates to deal with the case $p\in[3,5)$. The main
part of the proof of Theorem \ref{th1} is written with arguments coming from the
dynamical study of PDEs.
They are simpler than the corresponding ones of control theory, but far less
accurate since they do not give any estimation for the time of observability.
Anyway, such accuracy is not important here
since we use the unique continuation property for \eqref{eq-iii} with
$T=+\infty$. We briefly recall in 
Section \ref{sect-control} how these propagation of compactness and regularity
properties could have been proved with some arguments more usual in the
control theory. 

Moreover, we give two applications of our results in both contexts of control
theory and dynamical systems. First, as it is usual in control theory, some
results of stabilization can be coupled with local control theorems to provide
global controllability in large time.
\begin{theorem} \label{thmcontrol}
Assume that $f$ satisfies the conditions of Theorem \ref{th1}
or belongs to the generic set $\Gg$ defined by Theorem \ref{th2}.   
Let $R_0>0$ and $\omega$ satisfying the geometric control condition. Then, there
exists $T>0$ such that for any $(u_0,u_1)$ and $(\tilde{u}_0 ,\tilde{u}_1)$ in
$H^1_0(\Omega)\times L^2(\Omega)$, with 
\bna
\nor{(u_0,u_1)}{H^1\times L^2} \leq R_0; &\quad
&\nor{(\tilde{u}_0,\tilde{u}_1)}{H^1\times L^2}\leq R_0\\ 
\ena
there exists $g\in L^{\infty}([0,T],L^2(\Omega))$ supported in $[0,T]\times
\omega$ such that the unique strong solution of  
\begin{eqnarray*}
\left\lbrace
\begin{array}{rcl}
\Box u+\beta u+f(u)&=& g\quad \textnormal{on}\quad [0,T]\times \Omega\\
(u(0),\partial_t u(0))&=&(u_0,u_1) .
\end{array}
\right.
\end{eqnarray*}
satisfies $(u(T),\partial_t u(T))=(\tilde{u}_0,\tilde{u}_1)$.
\end{theorem}

The second application of our results concerns the existence of a compact global attractor. 
A compact global attractor is a compact set, which is invariant by the flow of
the PDE and which attracts the bounded sets. The existence of such an attractor
is an important dynamical property because it roughly says that the dynamics
of the PDE may be reduced to dynamics on a compact set, which is often
finite-dimensional. See for example \cite{Hale-book} and \cite{Raugel} for a
review on this concept. Theorems \ref{th1} and \ref{th2} show that $\{0\}$ is a
global attractor for the damped wave equation \eqref{eq}. Of course, it is
possible to obtain a more complex attractor by considering an equation of the
type
\begin{equation}\label{eq-attrac}
\left\{\begin{array}{ll}
\partial^2_{tt} u +\gamma(x) \partial_t u=\Delta
u-\beta u - f(x,u)~~&(x,t)\in\Omega\times\RR_+~,\\
u(x,t)=0&(x,t)\in\partial\Omega\times \RR_+\\
(u,\partial_t u)=(u_0,u_1)\in H^1_0\times L^2&
\end{array}\right.
\end{equation}
where $f\in C^{\infty}(\Omega\times\RR,\RR)$ is real analytic with respect to
$u$ and satisfies the following properties. There 
exist $C>0$, $p\in [1,5)$ and $R>0$ such that for all $(x,u)\in\Omega\times
\RR$, 
\begin{equation}\label{hyp-f2}
|f(x,u)|\leq C(1+|u|)^p~,~~|f'_x(x,u)|\leq C(1+|u|)^p~,~~|f'_u(x,u)|\leq
C(1+|u|)^{p-1}
\end{equation}
\begin{equation}\label{hyp-f3}
\left(x\not\in B(x_0,R) \text{ or }|u|\geq R\right)~~\Longrightarrow~~ f(x,u)u\geq
0~.
\end{equation}
where $x_0$ denotes a fixed point of the manifold.
\begin{theorem}\label{th-attrac}
Assume that $f$ is as above. Then, the dynamical
system generated by \eqref{eq-attrac} in $H^1_0(\Omega)\times L^2(\Omega)$ is
gradient and admits a compact global attractor $\Ac$. 
\end{theorem}
Of course, we would get the same result for $f$ in a generic set
similar to the one of Theorem \ref{th2}.

\bigskip

We begin this paper by setting our main notations and recalling the basic
properties of Equation \eqref{eq} in Section \ref{sect-basic}. 
We recall the unique continuation property of Robbiano and Zuily in Section
\ref{sect-UCP},
whereas Sections \ref{sect-asymp-comp} and \ref{sect-HR} are concerned by the
asymptotic compactness and the asymptotic smoothing effect of the damped wave
equation. The proofs of our main results, Theorem \ref{th1} and \ref{th2}, are
given in Sections \ref{sect-th1} and \ref{sect-th2} respectively. An
alternative proof, using more usual arguments from control theory, is sketched
in Section \ref{sect-control}. Finally, Theorems \ref{thmcontrol} and
\ref{th-attrac} are discussed in Section \ref{sect-disc}.

\medskip

{\bf \noindent Acknowledgements:} the second author was financed by the ERC
grant {\it GeCoMethods} during part of the redaction of this article. Moreover,
both authors benefited 
from the fruitful atmosphere of the conference Partial differential equations,
optimal design and numerics in Benasque. We also would like to thank Mathieu
L\'eautaud for his remarks about the optimality of Hypothesis (GCC) in the
nonlinear context and Genevi\`eve Raugel for her help for removing
a non-natural hypothesis of Theorem \ref{th-attrac}.

%%%%%%%%%%%%%%%%%%%%%%%%%%%%%%%%%%%%%%%%%%%%%%%%%%%%%%%%%%%%%%%%%%%%%%%%%%%%
%%%%%%%%%%%%%%%%%%%%%%%%%%%%%%%%%%%%%%%%%%%%%%%%%%%%%%%%%%%%%%%%%%%%%%%%%%%%

\bigskip

\section{Notations and basic properties of the damped wave
equation}\label{sect-basic}

In this paper, we use the following notations:
$$U=(u,u_t)~,~~F=(0,f)~,~~ A=\left(\begin{array}{cc}
0& Id\\ \Delta-\beta &-\gamma \end{array}\right)~.$$ 
In this setting, \eqref{eq} becomes
$$\partial_t U(t)=AU(t)+F(U)~.$$
We set $X=H^1_0(\Omega)\times L^2(\Omega)$ and for $s\in[0,1]$, $X^s$
denotes the space $D((-\Delta+\beta)^{(s+1)/2})\times D((-\Delta+\beta)^{s/2})=(H^{1+s}(\Omega)\cap
H^1_0(\Omega))\times
H^s_0(\Omega)$. Notice that $X^0=X$ and $X^1=D(A)$ (even if $\gamma$ is only in $L^{\infty}$).

We recall that $E$ denotes the energy defined by \eqref{def-E}. We also
emphasize that \eqref{hyp-f} and the invertibility of $\Delta-\beta$ 
implies that a set is bounded in $X$ if and only if its energy $E$ is bounded.
Moreover, for all $E_0\geq 0$, there exists $C>0$ such that
\begin{equation}\label{equiv-norm}
\forall (u,v)\in X,~E(u,v) \leq
E_0~\Longrightarrow~\frac{1}{C}\left\|(u,v)\right\|^2_X\leq
E(u,v)\leq  C\left\|(u,v)\right\|^2_X
\end{equation}
To simplify some statements in the proofs, we assume without loss of
generality that $3<p<5$. It will avoid some meaningless statements with negative
Lebesgue exponents since $p=3$ is the exponent where Strichartz  estimates are
no more necessary and can be replaced by Sobolev embeddings.

We recall that $\Omega$ is endowed with a metric $g$. We denote by $d$ the
distance on $\Omega$ defined by 
$$d(x,y)=\inf\left\{l(c)\left|c\in \Cc^{\infty}([0,1],\Omega) \text{ with
} c(0)=x \text{ and }c(1)=y\right. \right\}$$ 
where $l(c)$ is the length of the path $c$ according to the metric $g$. A ball
$B(x,R)$ in $\Omega$ is naturally defined by
$$B(x,R)=\{y\in\Omega, d(x,y)<R\}~.$$
For instance, if $\Omega=\R^3 \setminus B_{\RR^3}(0,1)$, the distance
between $(0,0,1)$ and $(0,0,-1)$ is $\pi$ (and not $2$) and the ball
$B((0,0,1),\pi)$ has nothing to do with the classical ball
$B_{\RR^3}((0,0,1),\pi)$ of $\RR^3$.

\subsection{Cauchy problem}

The global existence and uniqueness of solutions of the subcritical wave 
equation \eqref{eq} with $\gamma\equiv 0$ has been studied by Ginibre and Velo
in \cite{Gi-Ve-1} and \cite{Gi-Ve-2}. Their method also applies for $\gamma\neq
0$ since this term is linear and well defined in the energy space $X$.
Moreover, their argument to prove uniqueness also yields the continuity of the
solutions with respect to the initial data. 

The central argument is the use of Strichartz estimates.
\begin{theorem}[Strichartz estimates]$~$\\
Let $T>0$ and $(q,r)$ satisfying 
\begin{equation}\label{eq-strichartz}
\frac{1}{q}+\frac{3}{r}=\frac{1}{2},\quad q\in [7/2,+\infty].
\end{equation}
There exists $C=C(T,q)>0$ such that for every $G\in 
L^1([0,T],L^2(\Omega))$ and every $(u_0,u_1)\in X$, the solution $u$ of
\bneq
\Box u+\gamma(x)\partial_t u&=&G(t)\\
(u,\partial_t u)(0)&=&(u_0,u_1)
\eneq
satisfies the estimate
$$
\nor{u}{L^{q}([0,T],L^{r}(\Omega))}\leq
C\left(\nor{u_0}{H^1(\Omega)}+\nor{u_1}{L^2(\Omega)}+\nor{G}{L^{1}([0,T],L^{2}
(\Omega))}\right).
$$
\end{theorem}

The result was stated in the Euclidean space $\R^3$ by Strichartz
\cite{Strichartz} and Ginibre and Velo with $q\in (2,+\infty]$. Kapitanski
extended the result to variable coefficients in \cite{KapitanskiStrichartzvar}.
On a bounded domain, the first estimates were proved by Burq, Lebeau and
Planchon
\cite{BLP} for $q\in [5,+\infty]$ and extended to a larger range by Blair, Smith
and Sogge in \cite{StrichartzSoggewave}. Note that, thanks to the
counterexamples of Ivanovici \cite{counterIvanov}, we know that we cannot
expect
some Strichartz estimates in the full range of exponents in the presence of
boundaries.

From these results, we deduce the estimates for the damped wave equation by
absorption for $T$ small enough. We can iterate the operation in a uniform
number of steps. Actually, for the purpose of the semilinear wave equation, it
is sufficient to consider the Strichartz estimate
$L^{\frac{2p}{p-3}}([0,T],L^{2p}(\Omega))$ which gives $u^p\in
L^{\frac{2}{p-3}}([0,T],L^{2}(\Omega))\subset L^{1}([0,T],L^{2}(\Omega))$
because $1<\frac{2}{p-3}<+\infty$.

\begin{theorem}[Cauchy problem]\label{th-Cauchy}$~$\\
Let $f$ satisfies \eqref{hyp-f}. Then, for any $(u_0,u_1)\in
X=H^1_0(\Omega)\times L^2(\Omega)$ there exists a unique solution $u(t)$ of
the subcritical damped wave equation \eqref{eq}. Moreover, this solution is
defined for all $t\in \R$ and its energy $E(u(t))$ is non-increasing in time. 

For any $E_0\geq 0$, $T\geq 0$ and $(q,r)$ satisfying \eqref{eq-strichartz},
there exists a constant $C$ such that, if $u$ is a solution of \eqref{eq} with
$E(u(0))\leq E_0$, then 
$$\nor{u}{L^{q}([0,T],L^{r}(\Omega))}\leq
C\left(\nor{u_0}{H^1(\Omega)}+\nor{u_1}{L^2(\Omega)}\right)~.$$

In addition, for any $E_0\geq 0$ and $T\geq 0$, there exists a constant $C$
such that, if $u$ and $\tilde u$ are two solutions of \eqref{eq} with
$E(u(0))\leq E_0$ and $E(\tilde u(0))\leq E_0$, then
$$\sup_{t\in[-T,T]} \|(u,\partial_t u)(t)-(\tilde u,\partial_t \tilde
u)(t)\|_X\leq C \|(u,\partial_t u)(0)-(\tilde u,\partial_t \tilde
u)(0)\|_X~.$$
\end{theorem}
\begin{demo}
The existence and uniqueness for small times is a consequence of Strichartz
estimates and of the subcriticality of the nonlinearity, see \cite{Gi-Ve-2}. The
solution can be globalized backward and forward in time thanks to the energy
estimates (\ref{energydecay}) for smooth solutions. Indeed,
$$
E(t)\leq E(s)+C\int_t^s E(\tau)\, d\tau
$$
and thus, Gronwall inequality for $t\leq s$ and the decay of energy for $t\geq
s$ show that the energy does not blow up in finite time. This allows to extend
the solution for all times since the energy controls the norm of the space $X$
by (\ref{equiv-norm}). 

For the uniform continuity estimate, we notice that $w=u-\tilde{u}$
is solution of
\bneq
\Box w+\beta w+\gamma(x)\partial_t w&=&-wg(u,\tilde{u}) \\
(w,\partial_t w)(0)&=&(u,\partial_t u)(0)-(\tilde u,\partial_t \tilde
u)(0)
\eneq
where $g(s,\tilde{s})=\int_0^1 f'(s+\tau(\tilde{s}-s))\, d\tau$ fulfills
$|g(s,\tilde{s})|\leq C(1+|s|^{p-1}+|\tilde{s}|^{p-1})$. Let $q=\frac{2p}{p-3}$,
Strichartz and H\"older estimates give
\bna
\|(w,\partial_t w)(t)\|_{L^\infty([0,T],X)\cap L^q([0,T],L^{2p})}&\leq&
C\|(w,\partial_t w)(0)\|_X+C \nor{wg(u,\tilde{u})}{L^1([0,T],L^2)}\\
&\leq&C\|(w,\partial_t w)(0)\|_X+CT
\nor{w}{L^\infty([0,T],L^2)}\\ 
&&+T^{\theta}\nor{w}{L^q([0,T],L^{2p})}\left(\nor{u
}{L^q([0,T],L^{2p})}^{p-1}+\nor{\tilde{u}}{L^q([0,T],L^{2p})}^{p-1}\right)
\ena
with $\theta=\frac{5-p}{2}>0$. We get the expected result for
$T$ small enough by absorption since we already know a uniform bound
(depending on $E_0$) for the Strichartz norms of $u$ and $\tilde u$. Then, we
iterate the operation to
get the result for large $T$.\enp

\subsection{Exponential decay of the linear semigroup}

In this paper, we will strongly use the exponential decay for the linear
semigroup in the case where
$\gamma$ may vanish but satisfies the geometric assumptions of this
paper. In this case, \eqref{hyp-gamma} enables to control the decay of energy
outside a large ball and the geometric control condition
(GCC) enables to control the energy trapped in this ball.
\begin{prop}\label{prop-dec-expo}
Assume that $\gamma\in L^\infty(\Omega)$ satisfies \eqref{hyp-gamma} and (GCC). 
There exist two positive constants $C$ and $\lambda$ such that
$$\forall s\in[0,1]~,~\forall t\geq 0~,~~|||e^{At}|||_{\Lc(X^s)}\leq
Ce^{-\lambda t}~.$$
\end{prop}

The exponential decay of the damped wave equation under the geometric control
condition is well known since the works of
Rauch and Taylor on a compact manifold \cite{RT} and Bardos, Lebeau and Rauch
\cite{BLR,BLRinterne} on a bounded domain. Yet, we did not find any reference
for unbounded domain (\cite{AK} and \cite{Khen} concern unbounded domains but
local energy only). It is noteworthy that the decay of the
linear semigroup in unbounded domains seems not to have been extensively studied
for the moment.

We give a proof of Proposition \ref{prop-dec-expo} using microlocal defect
measure as done in Lebeau \cite{Lebeauamorties} or Burq \cite{Burqreg} (see also
\cite{BurqGerardCNS} for the proof of the necessity). The only difference with
respect to these results is that the manifold that we consider may be unbounded.
Since microlocal defect measure only reflects the local propagation, we thus
have to use the property of equipartition of the energy to deal with the
energy at infinity and to show a propagation of compactness (see \cite{DLZ}
for the flat case). 
\begin{lemma}\label{lemmepropag}
Let $T>L$ where $L$ is given by (GCC).
Assume that $(U_{n,0})\subset X$ is a bounded sequence, which weakly converges
to $0$ and assume that $U_n(t)=(u_n(t),\partial_t u_n(t))=e^{At}U_{n,0}$
satisfies
\bnan
\label{obsabsurd}
\int_0^T \int_{\Omega}\gamma(x)|\partial_t u_n|^2\rightarrow 0~. 
\enan
Then, $(U_{n,0})$ converges to $0$ strongly in $X$.
\end{lemma}
\bnp ~Let $\mu$ be a microlocal defect measure associated to $(u_n)$ (see
\cite{defectmeasure}, \cite{tartarh-measure} or \cite{BurqBourbaki} for the
definition). Note that (\ref{obsabsurd}) implies that $\mu $ can also be
associated to the solution of the wave equation without damping, so the weak
regularity of $\gamma$ is not problematic for the propagation and we get that
$\mu$ is 
concentrated on $\left\{\tau^2-|\xi|_x^2=0\right\}$ where $(\tau,\xi)$ are the
dual variables of $(t,x)$. Moreover,  (\ref{obsabsurd}) implies that
$\gamma\tau^2\mu=0$ and so $\mu\equiv 0$ on $S^*(]0,T[\times \omega)$. Then,
using the propagation of the measure along the generalized bicharacteristic flow
of Melrose-Sj\"ostrand and the geometric control condition satisfied by
$\omega$, we obtain $\mu\equiv 0$ everywhere. We do not give more details about
propagation of microlocal defect measure and refer to the Appendix
\cite{Lebeauamorties} or Section 3 of \cite{BurqBourbaki} (see also
\cite{GerardLeichtnam} for some close propagation results in a different
context). Since $\mu \equiv 0$, we know that $U_n\rightarrow 0$ on $H^1\times
L^2(]0,T[\times B(x_0,R))$ for every $R>0$. 

To finish the proof, we need the classical equipartition of the energy to get
the convergence to $0$ in the whole manifold $\Omega$. Since $\gamma$ is
uniformly positive outside a ball $B(x_0,R)$, \eqref{obsabsurd} and the
previous arguments imply that $\partial_t u_n \rightarrow 0$ in $L^2([0,T]
\times \Omega)$. Let $\varphi \in C^{\infty}_0(]0,T[)$ with $\varphi\geq 0$ and
$\varphi(t)=1$ for $t\in [\e,T-\e]$. We multiply the equation by $\varphi(t)u_n$
and we obtain
\bna
0&=&-\iint_{[0,T]\times \Omega} \varphi(t)|\partial_t u_n|^2-\iint_{[0,T]\times
\Omega} \varphi'(t) \partial_t u_n u_n +\iint_{[0,T]\times \Omega}
\varphi(t)|\nabla u_n|^2\\
&&+\iint_{[0,T]\times \Omega}  \varphi(t)\beta|u_n|^2+\iint_{[0,T]\times
\Omega}\varphi(t) \gamma(x)\partial_t u_n u_n. 
\ena
The $L^2-$norm of $u_n(t)$ is bounded, while $\partial_t u_n \rightarrow 0$ in
$L^2([0,T] \times \Omega)$, so the first, second and fifth terms converge to
zero. Then, the above equation yields
\bna
\iint_{[0,T]\times \Omega} \varphi(t)\left(\beta|u_n|^2+|\nabla
u_n|^2\right)\longrightarrow 0.
\ena
Finally, notice that the energy identity
$\nor{U_{n,0}}{X}^2=\nor{U_{n}(t)}{X}^2+\int_0^T
\int_{\Omega}\gamma(x)|\partial_t u_n|^2$ shows that 
$$\iint_{[0,T]\times \Omega} \varphi(t)\left(\beta|u_n|^2+|\nabla u_n|^2\right)
\sim \nor{U_{n,0}}{X}^2 \int_0^T \varphi(t)$$
and thus that $\nor{U_{n,0}}{X}$ goes to zero. 
\enp

\noindent \emph{\textbf{Proof of Proposition  \ref{prop-dec-expo}:}} Once Lemma
\ref{lemmepropag} is established, the proof follows the arguments of the
classical case, where $\Omega$ is bounded. We briefly recall them.

We first treat the case $s=0$. As in Proposition \ref{critere-decay}, the
exponential decay of the energy is equivalent to the observability estimate,
that is the existence of $C>0$ and $T>0$ such that, for any trajectory
$U(t)=e^{At}U_0$ in $X$,
\bna
\label{eq-critere-decaylin}
\int_0^T \int_{\Omega}\gamma(x)|\partial_t u|^2 \geq C\nor{U(0)}{X}^2~.
\ena

We argue by contradiction: assume that
\eqref{eq-critere-decaylin} does not hold for any positive $T$ and $C$.
Then, there exists a sequence of initial data $U_n(0)$ with
$\|U_n(0)\|_X=1$ and such that 
$$\int_0^n \int_\Omega \gamma(x)|\partial_t
u_n(t,x)|^2 dtdx \xrightarrow[~~n\longrightarrow +\infty~~]{} 0~,$$
where $(u_n,\partial_t u_n)(t)=U_n(t)=e^{At}U_n(0)$. 
Let $\tilde U_n=U_n(n/2+\cdot)$. We have 
$$\int_{-n/2}^{n/2} \int_\Omega
\gamma(x)|\partial_t \tilde u_n(t,x)|^2 dtdx \xrightarrow[~n\rightarrow
\infty~]{} 0~,$$
and, for any $t\in[-n/2,n/2]$,
$$\|\tilde U_n(t)\|_X^2 = \|\tilde U_n(-n/2)\|_X^2- \int_{-n/2}^{t}
\int_\Omega \gamma(x)|\partial_t \tilde u_n(s,x)|^2 dsdx
\xrightarrow[~n\rightarrow
\infty~]{} 1~.$$
Therefore, we can assume that $U_n(0)$ converges to $U_\infty(0)\in X$,
weakly in $X$. Moreover, for any $T>0$, $U_n(t)$ and $\partial_t U_n(t)$ are
bounded in $L^\infty([-T,T],X)$ and $L^\infty([-T,T],L^2(\Omega)\times
H^{-1}(\Omega))$ respectively. Thus, using Ascoli's Theorem, we may also assume
that $U_n(t)$ strongly converges to $U_\infty(t)$ in
$L^\infty([-T,T], L^2(K)\times H^{-1}(K))$ where $K$ is
any compact of $\Omega$. Hence, $(u_\infty,\partial_t
u_\infty)(t)=U_\infty(t)=e^{At}U_\infty(0)$ is a solution of 
\bneqn
\label{UCPlin}
\Box u_{\infty}+\beta u_{\infty}&=&0\textnormal{ on } \RR\times \Omega\\
\partial_t u_{\infty}&= &0 \textnormal{ on } \RR\times \omega.
\eneqn
in $L^2\times H^{-1}$. Since $U_\infty(0)\in X$ belongs to $X$, we deduce that,
in fact, $U_\infty(t)$ solves \eqref{UCPlin} in $X$. 

To finish the proof of Proposition \ref{prop-dec-expo}, we have to show that
$U_\infty\equiv 0$. Indeed, applying Lemma \ref{lemmepropag}, we would get that
$U_n$ converges strongly to $0$, which contradicts the hypothesis
$\|U_n(0)\|_X=1$. Note that $U_\infty\equiv 0$  is a direct consequence of a
unique continuation property as Corollary \ref{coro-RZ}. However, 
Corollary \ref{coro-RZ} requires $\Omega$ to be smooth, whereas Proposition
\ref{prop-dec-expo} could be more general. Therefore, we recall another
classical argument to show that $U_\infty\equiv 0$.

Denote $N$ the set of function $U_{\infty}(0)\in X$ satisfying \eqref{UCPlin},
which is obviously a linear subspace of $X$. We will prove that
$N=\left\{0\right\}$. Since $\gamma(x)|\partial_t u_\infty|^2\equiv
0$ for functions $u_\infty$ in $N$ and since $N$ is a closed subspace, Lemma
\ref{lemmepropag} shows that any weakly convergent subsequence of $N$ is in fact
strongly convergent. By Riesz Theorem, $N$ is therefore finite dimensional. For
any $t\in \R$, $e^{tA}$ applies $N$ into itself and thus $A_{|N}$ is a bounded
linear operator. Assume that $N\neq \{0\}$, then $A_{|N}$ admits an
eigenvalue $\lambda$ with eigenvector $Y=(y_0,y_1)\in N$. This means that
$y_1=\lambda y_0 $ and that $(\Delta-\beta)y_0=\lambda^2 y_0$. 
Moreover, we know that $y_1=0$ on $\omega$ and so, if $\lambda\neq 0$,
that $y_0=0$ on $\omega$. This implies $y_0\equiv 0$ by the unique continuation
property of elliptic operators. Finally, if $\lambda=0$, we have
$(\Delta-\beta)y_0=0$ and $y_0=0$, because, by assumption, $\Delta-\beta$ is a
negative definite operator. 

So we have proved $N=\left\{0\right\}$ and therefore $U_{\infty}=0$, that is
$\tilde U_n(0)$ converges to $0$
weakly in $X$. We can then apply Lemma \ref{lemmepropag} on any interval
$[-n/2,-n/2+T]$ where $L$ is the time in the geometric control condition (GCC)
and obtain a contradiction to $\|U_n(0)\|_X=1$. 

Let us now consider the cases $s\in (0,1]$. The basic semigroup properties (see
\cite{Pazy}) shows that, if $U\in
X^1=D(A)$, then $e^{At}U$ belongs to $D(A)$ and 
\bna
\|e^{At}U\|_{X^1}&=&\|Ae^{At}U\|_X+\|e^{At}U\|_X=\|e^{At}AU\|_X+\|e^{At}U\|_X\\
&\leq &Ce^{-\lambda t} \left(\|AU\|_X+\|U\|_X\right)
= Ce^{-\lambda t} \|U\|_{D(A)}~.
\ena
This shows Proposition \ref{prop-dec-expo} for $s=1$. Notice that we do not
have to require any regularity for $\gamma$ to obtain this result. Then,  
Proposition \ref{prop-dec-expo} for $s\in (0,1)$ follows by interpolating
between the cases $s=0$ and $s=1$ (see \cite{Tartar}).
{\hfill$\square$\\
\vspace{0.4cm}}

\subsection{First nonlinear exponential decay properties}\label{sect-exp-dec}

Theorem \ref{th-Cauchy} shows that the energy $E$ is non-increasing along the
solutions of \eqref{eq}. The purpose of this paper is to obtain the exponential
decay of this energy in the sense of property (ED) stated above. We first recall
the well-known criterion for exponential decay.
\begin{prop}\label{critere-decay}
The exponential decay property (ED) holds if and only if,
there exist $T$ and $C$ such that 
\begin{equation}\label{eq-critere-decay}
E(u(0))\leq C(E(u(0))-E(u(T)))= C \int\int_{[0,T]\times\Omega}
\gamma(x)|\partial_t u(x,t)|^2 dtdx
\end{equation}
for all solutions $u$ of \eqref{eq} with $E(u(0))\leq E_0$.
\end{prop}
\begin{demo}
If (ED) holds then obviously \eqref{eq-critere-decay} holds for $T$ large
enough since $E(u(0))-E(u(T))\geq (1-Ke^{-\lambda T})E(u(0))$. Conversely, if
\eqref{eq-critere-decay} holds, using $E(u(T))\leq E(u(0))$,
we get $E(u(T))\leq C/(C+1)E(u(0))$ and thus $E(u(kT))\leq
(C/(C+1))^k E(u(0))$. Using again the decay of the energy to fill the gaps
$t\in(kT,(k+1)T)$, this shows that (ED) holds. 
\end{demo}

First, we prove exponential decay in the case of positive damping, which will be
helpful to study
what happens outside a large ball since \eqref{hyp-gamma} is assumed in the
whole paper. Note that the fact that $-\Delta+\beta$ is positive is necessary to
avoid for instance the constant undamped solutions.
\begin{prop}\label{prop-positive-damping}
Assume that $\omega=\Omega$, that is that $\gamma(x)\geq \alpha>0$
everywhere. Then (ED) holds.
\end{prop}
\begin{demo}
We recall here the classical proof. We introduce a modified energy
$$\tilde E(u)=\int_\Omega \frac12 (|\partial_t u|^2+|\grad u|^2+\beta
|u|^2)+V(u) + \varepsilon u\partial_tu$$
with $\varepsilon>0$. Since $\int_\Omega |\grad u|^2+\beta |u|^2$ controls
$\|u\|_{L^2}^2$, $\tilde E$ is equivalent to $E$ for $\varepsilon$ small
enough and it is sufficient to obtain the exponential decay of the auxiliary
energy $\tilde E$. Using $\gamma\geq \alpha>0$ and $uf(u)\geq 0$, a direct
computation shows for $\varepsilon$ small enough that
\begin{align*}
\tilde E(u(T))-\tilde E(u(0))&=\int_0^T\int_\Omega -\gamma(x)|\partial_t u|^2 +
\varepsilon |\partial_t u|^2+\varepsilon \gamma(x)u\partial_t u -\varepsilon
(|\grad u|^2+\beta|u|^2)-\varepsilon uf(u)\\
&\leq - C \int_0^T \|(u,\partial_t u)\|_{H^1\times L^2}^2\\
&\leq -C \int_0^T \tilde E(t)dt\\
&\leq -CT \tilde E(T)~,
\end{align*}
where $C$ is some positive constant, not necessarily the same from line to line.
Thus, $\tilde E(u(0))-\tilde E(u(T))\geq CT \tilde E(u(T))$ with $CT>0$
and therefore $\tilde E(u(0))\geq \mu \tilde E(u(T))$ with $\mu>1$. As in
the proof of Proposition \ref{prop-dec-expo}, this last property implies the
exponential decay of $\tilde E$ and thus the one of $E$.
\end{demo}

%%%%%%%%%%%%%%%%%%%%%%%%%%%%%%%%%%%%%%%%%%%%%%%%%%%%%%%%%%%%%%%%%%%%%%%%%%%%
%%%%%%%%%%%%%%%%%%%%%%%%%%%%%%%%%%%%%%%%%%%%%%%%%%%%%%%%%%%%%%%%%%%%%%%%%%%%

\bigskip      

\section{A unique continuation result for equations with partially
holomorphic coefficients}\label{sect-UCP}

Comparatively to previous articles on the stabilization of the damped wave
equations as \cite{DLZ}, one of the main novelties of this paper is the use of
a unique continuation theorem requiring partially analyticity of the
coefficients, but very weak geometrical assumptions as shown in
Corollary \ref{coro-RZ}. We use here the
result of Robbiano and Zuily in \cite{Rob-Zui}. This result has also been proved
independently by H\"ormander in \cite{Hormander_uniq_analytic} and has been
generalised by Tataru in \cite{Tataru_uniq_analyticJMPA}. Note that the idea
of using partial analyticity for unique continuation was introduced by Tataru
\cite{Tataru_uniq_analyticCPDE} but it requires some global analyticity
assumptions that are not fulfilled in our case. All these results use
very accurate microlocal analysis and hold in a much more general framework
than the one of the wave equation. However, for sake of simplicity, we restrict
the statement to this case.
\begin{theorem}\label{th-RZ}
{\bf Robbiano-Zuily, H\"ormander (1998)}\\
Let $d\geq 1$, $(x_0,t_0)\in\RR^d\times\RR$ and let $\Uc$ be a neighbourhood of 
$(x_0,t_0)$. Let $(A_{i,j}(x,t))_{i,j=1,\ldots,d}$, $b(x,t)$,
$(c_i(x,t))_{i=1,\ldots,d})$ and $d(x,t)$ be bounded coefficients in
$\Cc^\infty(\Uc,\RR)$. Let $v$ be a strong solution of 
\begin{equation}\label{eq-RZ}
\partial^2_{tt} v=\div(A(x,t)\grad v) + b(x,t)\partial_t v+c(x,t).\grad
v+d(x,t)v~~~~~(x,t)\in\Uc\subset\RR^d\times\RR~.
\end{equation}
Let $\varphi\in\Cc^2(\Uc,\RR)$ such that $\varphi(x_0,t_0)=0$ and $(\grad
\varphi,\partial_t\varphi)(x,t)\neq 0$ for all $(x,t)\in\Uc$. Assume that:\\
(i) the coefficients $A$, $b$, $c$ and $d$ are analytic in time,\\
(ii) $A(x_0,t_0)$ is a symmetric positive definite matrix,\\
(iii) the hypersurface $\{(x,t)\in\Uc~,~\varphi(x,t)=0\}$ is not characteristic
at $(x_0,t_0)$ that is that we have $|\partial_t\varphi(x_0,t_0)|^2\neq \langle
\grad
\varphi(x_0,t_0) |A(x_0,t_0) \grad\varphi (x_0,t_0)\rangle$\\
(iv) $v\equiv 0$ in $\{(x,t)\in\Uc~,~\varphi(x,t)\leq 0\}$.\\
Then, $v\equiv 0$ in a neighbourhood of $(x_0,t_0)$.
\end{theorem}
\begin{demo}
We only have to show that Theorem \ref{th-RZ} is a direct
translation of Theorem A of \cite{Rob-Zui} in the framework of the wave
equation. To use the notations of \cite{Rob-Zui}, we set
$x_a$ to be the time variable and $x_b$ the space variable and we set
$(x_0,t_0)=x^0=(x_b^0,x_a^0)$. Equation \eqref{eq-RZ} corresponds
to the differential operator 
$$P=\xi_a^2-\,^t\xi_b A(x_b,x_a) \xi_b -
b(x_b,x_a)\xi_a-c(x_b,x_a)\xi_b-d(x_b,x_a)$$
with principal symbol $p_2=\xi_a^2-\,^t\xi_b A(x_b,x_a) \xi_b $. 

All the statement of Theorem \ref{th-RZ} is an obvious
translation of Theorem A of \cite{Rob-Zui}, except maybe for the fact
that Hypothesis (iii) implies the hypothesis of pseudo-convexity of
\cite{Rob-Zui}. We compute $\left\{p_2,\varphi\right\}=2\xi_a\varphi'_a -2^t\xi_b A(x_a,x_b) \varphi'_b$. Let us set
$\zeta=(x_a^0, x_b^0, i\varphi'_a(x^0), \xi_b+i\varphi'_b(x^0))$, then
$\left\{p_2,\varphi\right\}(\zeta)=0$ if and only if
$$i(\varphi'_a(x^0))^2-i^t\varphi'_b(x^0))A(x^0)\varphi'_b(x^0)-
^t\xi_bA(x^0)\varphi'_b(x^0)=0~.$$
This is possible only if $(\varphi'_a(x^0))^2=\,^t\varphi'_b(x^0)A(x^0)\varphi'_b(x^0)$, that is if the
hypersurface
$\varphi=0$ is characteristic at $(x_0,t_0)$. Thus, if this hypersurface is not
characteristic, then the pseudo-convexity hypothesis of Theorem A of
\cite{Rob-Zui} holds. 
\end{demo}

The previous theorem allows to prove some unique continuation result with some
optimal time and geometric assumption. This allows to prove unique continuation
where the geometric condition is only, roughly speaking, that we do not
contradict the finite speed of propagation. 
\begin{corollary}\label{coro-RZ}
Let $T>0$ (or $T=+\infty$) and let $b$, $(c_i)_{i=1,2,3}$ and
$d$ be smooth coefficients in $\Cc^\infty(\Omega\times [0,T],\RR)$. Assume
moreover that $b$, $c$
and $d$ are analytic in time and that $v$ is a strong solution of
\begin{equation}\label{eq-RZ-2}
\partial^2_{tt} v=\Delta v + b(x,t)\partial_t v+c(x,t).\grad
v+d(x,t)v~~~~~(x,t)\in\Omega\times (-T,T)~.
\end{equation}
Let $\Oc$ be a non-empty open subset of $\Omega$ and assume that 
$v(x,t)=0$ in $\Oc\times (-T,T)$. Then
$v(x,0)\equiv 0$ in $\Oc_T=\{x_0\in\Omega~,~d(x_0,\Oc)<T\}$.\\
As consequences:\\
a) if $T=+\infty$, then $v\equiv 0$ everywhere,\\
b) if $v\equiv 0$ in $\Oc\times (-T,T)$ and $\overline \Oc_T=\Omega$, then
$v\equiv 0$ everywhere.
\end{corollary}
\begin{demo}
Since $\Omega$ is assumed to be connected, both consequences are obvious from
the first statement.  

Let $x_0$ be given such that $d(x_0,\Oc)<T$. There is a point $x_*\in\Oc$ linked
to $x_0$ by a smooth curve of length $l<T$, which stays away from the boundary.
We introduce a sequence of balls $B(x_0,r)$, \ldots, $B(x_K,r)$ with
$r\in(0,T/K)$, $x_{k-1}\in B(x_{k},r)$ and $x_K=x_*$, such that
$B(x_k,r)$ stays away from the boundary and is
small enough such that it is diffeomorphic to an open set of $\RR^3$ via the
exponential map. Note that such a sequence of balls exists because the smooth
curve linking $x_0$ to $x_K$ is compact and of length smaller than $T$. We also
notice that it is sufficient to prove Corollary \ref{coro-RZ} in each ball
$B(x_k,r)$. Indeed, this would enable us to apply Corollary \ref{coro-RZ} in
$B(x_K,r)\times (-T,T)$ to obtain that $v$ vanishes in a neighbourhood of
$x_{K-1}$ for $t\in (-T+r,T-r)$ and then to apply it recursively in
$B(x_{K-1},r)\times (-T+r,T-r)$, \ldots, $B(x_1,r)\times(-T+(K-1)r,T-(K-1)r)$ to
obtain that $v(x_0,0)=0$.

From now on, we assume that $x_0\in B(x_*,r)$ and that $v$ vanishes in a
neighbourhood $\Oc$ of $x_*$ for $t\in (-r,r)$. Since $d(x_0,x_*)<r$, we can
introduce a non-negative function $h\in\Cc^\infty([-r,r],\RR)$ such that
$h(0)>d(x_0,x_*)$, $h(\pm r)=0$ and $|h'(t)|< 1$ for all $t\in[-r,r]$. 
We set $\Uc=B(x_*,r)\times(-r,r)$ and for any $\lambda\in [0,1]$, we define 
$$\varphi_\lambda(x,t)=d(x,x_*)^2-\lambda h(t)^2~.$$
Since $r$ is assumed to be smaller than the radius of injectivity of the
exponential map, $\varphi_\lambda$ is a smooth well-defined function. 
We prove Corollary \ref{coro-RZ} by contradiction. Assume that $v(x_0,0)\neq 0$.
We denote by $V_{\lambda}$ the volume $\{(x,t)\in \Uc,~\varphi_\lambda(x,t)\leq
0\}$. We notice that $V_{\lambda_1}\subset V_{\lambda_2}$ if
$\lambda_1<\lambda_2$, that for small
$\lambda$, $V_\lambda$ is included in $\Oc\times (-r,r)$ where $v$
vanishes,
and that $V_1$ contains $(x_0,0)$ where $v$ does not vanish. Thus 
$$\lambda_0=\sup \{\lambda\in [0,1]~,~ \forall (x,t)\in V_\lambda,~ v(x,t)=
0\}~$$
is well defined and belongs to $(0,1)$. For $t$ close to $-r$ or $r$,
$h(t)$ is small and the section $\{x,~(x,t)\in V_{\lambda_0}\}$ of
$V_{\lambda_0}$ is contained in $\Oc$ where $v$ vanishes. Therefore, by
compactness, the hypersurface $S_{\lambda_0}=\partial V_{\lambda_0}$ must touch
the support of $v$ at some point $(x_1,t_1)\in\Uc$ (see Figure \ref{figRZ}).

\begin{figure}[htp]
\begin{center}
\resizebox{0.4\textwidth}{!}{\input{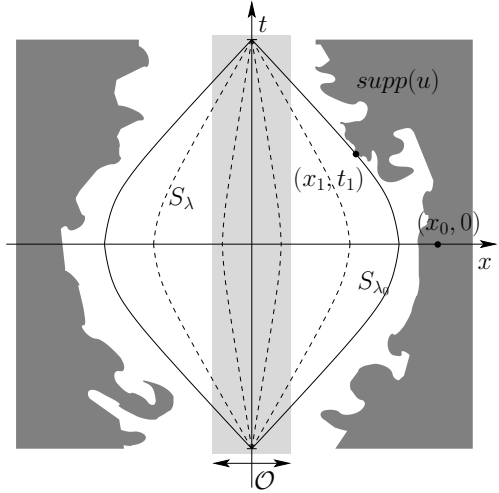}}
\end{center}
\caption{\it The proof of Corollary \ref{coro-RZ}} 
\label{figRZ}
\end{figure}

In local coordinates, $\Delta$ can be written $\div(A(x)\nabla
.)+c(x)\cdot \nabla $. Moreover, $ \langle \nabla
\varphi_\lambda | A \nabla \varphi_\lambda \rangle = |\nabla_g d(.,x_*)|_g^2=1$
where the index $g$ means that the gradient and norm are taken according to the
metric. Therefore,  
the hypersurface $S_{\lambda_0}$ is non-characteristic at $(x_1,t_1)$ in the
sense of Hypothesis (iii) of Theorem \ref{th-RZ} since $|\partial_t
\varphi_\lambda(x,t)|=|\lambda h'(t_1)|< 1$.
Thus, we can apply Theorem \ref{th-RZ} with $\varphi=\varphi_{\lambda_0}$ at the
point
$(x_1,t_1)$, mapping everything in the 3d-euclidean frame via the exponential
chart. We get that $v$ must vanish in a neighbourhood of $(x_1,t_1)$. This is
obviously a contradiction since $(x_1,t_1)$ has been taken in the support of
$v$.
\end{demo}

%%%%%%%%%%%%%%%%%%%%%%%%%%%%%%%%%%%%%%%%%%%%%%%%%%%%%%%%%%%%%%%%%%%%%%%%%%%%
%%%%%%%%%%%%%%%%%%%%%%%%%%%%%%%%%%%%%%%%%%%%%%%%%%%%%%%%%%%%%%%%%%%%%%%%%%%%

\bigskip      

\section{Asymptotic compactness}\label{sect-asymp-comp}

As soon as $t$ is positive, a solution $u(t)$ of a parabolic PDE becomes smooth
and stays in a compact set. The smoothing effect in finite time of course fails
for the damped wave equations. However, these PDEs admit in some sense a
smoothing effect in infinite time. This effect is called {\it asymptotic
compactness} if one is interested in extracting asymptotic subsequences as in
Proposition \ref{prop-asympt-reg}, or {\it asymptotic smoothness} if one uses
the regularity of globally bounded solutions as in Proposition
\ref{prop-reg-1}. For the reader interested in these notions, we refer for
example to \cite{Hale-book}. The proof of this asymptotic smoothing effect is
based on the variation of constant formula $U(t)=e^{At}U_0+\int_0^t
e^{A(t-s)}F(U(s))ds$ and two properties:\\ 
- the exponential decay of the linear group (Proposition \ref{prop-dec-expo}),
which implies that the linear part $e^{At}U_0$ asymptotically disappears,\\
- the regularity of the nonlinearity $F$ implying the compactness of the
nonlinear term $\int_0^t e^{A(t-s)}F(U(s))ds$ (Corollary \ref{coro-DLZ} below).
Note that the subcriticality of $f$ is the key point of this property and that
our arguments cannot be extended as they stand to the critical case $p=5$.

The purpose of this section is to prove some compactness and regularity results
about undamped solutions as (\ref{eq-iii}). Note that these results could also
have been obtained with a more ``control theoretic'' proof (see section
\ref{sect-control} for a sketch of the alternative proof) based on propagation
results or observability estimates. Here, we have chosen to give a different one
using asymptotic regularization, more usual in  dynamical system.  The spirit of
the proof remains quite similar: prove that the nonlinearity is more regular
than it seems {\it a priori} and use some properties of the damped linear
equation. 

\subsection{Regularity of the nonlinearity}

Since $f$ is subcritical, it is shown in \cite{DLZ} that the nonlinear term of
\ref{eq} yields a gain of smoothness.
\begin{theorem}\label{th-DLZ}
{\bf Dehman, Lebeau and Zuazua (2003)}\\
Let $\chi\in \Cc^\infty_0(\RR^3,\RR)$, $R>0$ and $T>0$. Let $s\in [0,1)$ and
let $\varepsilon=\min(1-s,(5-p)/2,(17-3p)/14)>0$ with $p$ and $f$ as in
\eqref{hyp-f}. There exist $(q,r)$ satisfying \eqref{eq-strichartz} and $C>0$
such that the following property holds.
If $v\in L^\infty([0,T],H^{1+s}(\RR^3))$ is a
function with finite Strichartz norms $\nor{v}{L^q([0,T],L^r(\RR^3))} \leq R$,
then
$\chi(x)f(v)\in L^1([0,T],H^{s+\varepsilon}(\RR^3))$ and moreover
$$\|\chi(x)f(v)\|_{L^1([0,T],H^{s+\varepsilon}(\RR^3))}\leq C
\|v\|_{L^\infty([0,T],H^{1+s}(\RR^3))}~.$$ 
The constant $C$ depends only on $\chi$, $s$, $T$, $(q,r)$, $R$ and the
constant in Estimate \eqref{hyp-f}.
\end{theorem}

Theorem \ref{th-DLZ} is a copy of Theorem 8 of \cite{DLZ},
except for two points.

First, we would like to apply the result to a solution $v$ of the damped
wave equation on a manifold possibly with boundaries, where not all Strichartz
exponents are available. This leads to the constraint 
$q\geq 7/2$ for the Strichartz exponents $(q,r)$ of \eqref{eq-strichartz} (see
Theorem \ref{th-Cauchy}). In the proof of Theorem 8 of \cite{DLZ}, the useful
Strichartz estimate corresponds to $r=\frac{3(p-1)}{1-\e}$ and
$q=\frac{2(p-1)}{p-3+2\e}$ and it is required that $q\geq p-1$, which yields
$\e\leq (5-p)/2$. In this paper, we require also that $q\geq 7/2$ which yields
in addition $\e \leq (17-3p)/14$. Notice that $p<5$ and thus both bounds are
positive.

The second difference is that, in \cite{DLZ}, the function $f$ is assumed to be
of class $\Cc^3$ and to satisfy
\begin{equation}\label{eq-DLZ}
|f''(u)|\leq C(1+|u|)^{p-2} ~~\text{ and }~~|f^{(3)}(u)|\leq C(1+|u|)^{p-3}
\end{equation}
in addition of \eqref{hyp-f}. Since Theorem \ref{th-DLZ} concerns the
$L^1(H^{s'})-$norm of $\chi(x)f(v)$ for $s'=s+\varepsilon\leq 1$, we can
omit Assumption \eqref{eq-DLZ}. Actually, we make the assumption $\varepsilon\leq 1-s$ which is not present in \cite{DLZ} and a careful study of their proof shows that \eqref{hyp-f} is not necessary under that assumption. 

Indeed, let $\tilde f(u)=\th^3(u)|u|^p$.
The function $\tilde f$ is of class $\Cc^3$ and satisfies \eqref{hyp-f} and
\eqref{eq-DLZ}. Hence, Theorem 8 of \cite{DLZ} can be applied to $\tilde f$ and
we can bound the $L^1(H^{s'})-$norm of $\tilde f$ as in Theorem \ref{th-DLZ}.
On the other hand, we notice that $|\tilde f(u)|\underset{\pm\infty}{\sim}
|u|^p$, $\tilde f'(u)\underset{\pm\infty}{\sim}
p|u|^{p-1}$ and $\tilde f'(u)\geq 0$. Therefore,
since $f$ satisfies \eqref{hyp-f}, there exists $C>0$ such that $|f(u)|\leq C
(1+ |\tilde f(u)|)$ and $|f'(u)|\leq C (1+\tilde f'(u))$. Thus, if say $v>u$ for
fixing the notations, 
\begin{align*}
|f(v)-f(u)|&\leq (v-u) \int_0^1 |f'(u+\tau(v-u))|d\tau\\
&\leq C(v-u)+C(v-u)\int_0^1 \tilde f'(u+\tau(v-u)) d\tau = C
(v-u)+C(\tilde f(v)-\tilde f(u))\\
&\leq C |v-u|+C|\tilde f(v)-\tilde f(u)| ~.
\end{align*}
For $0<s<1$, using the above inequalities and the definition of the
$H^{s'}-$norm as
$$\|\chi f(u)\|_{H^{s'}}^2= \|\chi f(u)\|_{L^2}^2 + \iint_{\RR^6} 
\frac {|\chi(x)f(u(x))-\chi(y)f(u(y))|^2}{|x-y|^{2s'}} dxdy~,
$$                                                                       
we obtain 
$$\|\chi f(u)\|_{L^1(H^{s'})}\leq C \|u\|_{L^\infty(H^1)}+C\| \tilde{\chi} \tilde
f(u)\|_{L^1(H^{s'})} $$
where $\tilde{\chi}$ is another cut-off function with larger support.
Hence, for $0<s<1$, the conclusion of Theorem \ref{th-DLZ} holds not only for
$\tilde f$ but also for $f$. If $s'=1$, we just apply the chain rule and
the proof is easier.

Note that the above arguments show that the constant $C$ depends on $f$
through Estimate \ref{hyp-f} only. Notice in addition that, since $f$ is only
$\Cc^1$, we cannot expect $\chi f(v)$ to be more regular than $H^1$ and that is
why we also assume $\e\leq 1-s$.\\[3mm]

In this paper, we use a generalisation of Theorem 8 of \cite{DLZ} to
non-compact manifolds with boundaries.
\begin{corollary}\label{coro-DLZ}
Let $R>0$ and $T>0$. Let $s\in [0,1)$ and
let $\varepsilon=\min(1-s,(5-p)/2,(17-3p)/14)>0$ with $p$ as in
\eqref{hyp-f}. There exist $(q,r)$ satisfying \eqref{eq-strichartz} and $C>0$
such that the following property holds.
If $v\in L^\infty([0,T],H^{1+s}(\Omega)\cap H^1_0(\Omega))$ is a
function with finite Strichartz norms $\nor{v}{L^q([0,T],L^r(\Omega))} \leq R$,
then $f(v)\in L^1([0,T],H^{s+\varepsilon}_0(\Omega))$ and moreover
$$\|f(v)\|_{L^1([0,T],H^{s+\varepsilon}_0(\Omega))}\leq C
\|v\|_{L^\infty([0,T],H^{s+1}(\Omega)\cap H^1_0(\Omega))}~.$$
The constant $C$ depends only on $\Omega$, $(q,r)$, $R$ and the constant in
Estimate
\eqref{hyp-f}.
\end{corollary}
\begin{demo}
Since we assumed that $\Omega$ has a bounded geometry in the sense that
$\Omega$ is compact or a compact perturbation of a manifold with periodic
metric, $\Omega$ can be covered by a set of $\Cc^\infty-$charts
$\alpha_i:U_i\longmapsto \alpha_i(U_i)\subset\RR^3$ such that $\alpha_i(U_i)$ is
equal either to $B(0,1)$ or to $B_+(0,1)=\{x\in B(0,1), x_1>0\}$ and such that,
for any $s\geq 0$ the norm of a function $u\in H^s(\Omega)$ is equivalent to the
norm 
$$(\sum_{i\in\NN} \|u\circ\alpha_i^{-1}\|^2_{H^{s}(\alpha_i(U_i))})^{1/2}~.$$
Moreover, the Strichartz norm $L^q([0,T],L^r(\alpha_i(U_i))$ of
$v\circ\alpha_i^{-1}$ is uniformly controlled from above by the Strichartz
norm $L^q([0,T],L^r(U_i))$ of $v$ which is bounded by $R$.

Therefore, it is sufficient to prove that Corollary \ref{coro-DLZ} holds for 
$\Omega$ being either $B(0,1)$ or $B_+(0,1)$. Say that $\Omega=B_+(0,1)$, the
case $\Omega=B(0,1)$ being simpler. To apply Theorem \ref{th-DLZ}, we extend $v$
in a neighbourhood of $B_+(0,1)$ as follows. For $x\in B_+(0,2)$, we use the
radial coordinates $x=(r,\sigma)$ and we set 
$$\tilde v(x)=\tilde v(r,\sigma)= 5 v(1-r,\sigma)-20v(1-r/2,\sigma)+16
v(1-r/4,\sigma)~.$$
Then, for $x=(x_1,x_2,x_3)\in B_-(0,2)$, we set 
$$\tilde v(x)=5 v(-x_1,x_2,x_3)-20v(-x_1/2
,x_2,x_3)+16 v(-x_1/4,x_2,x_3)~.$$
Notice that $\tilde v$ is an extension of $v$ in $B(0,2)$, which preserves the
$\Cc^2-$regularity, and that the $H^s-$norm for $s\leq 2$ as well as the
Strichartz norms of $\tilde v$ are controlled by the corresponding norms of $v$.
Let $\chi\in\Cc^\infty_0(\RR^3)$ be a cut-off
function such that $\chi\equiv 1$ in $B_+(0,1)$ and $\chi\equiv 0$ outside
$B(0,2)$. Applying Theorem \ref{th-DLZ} to $\chi(x)f(\chi(x)\tilde v)$ yields a
control of $\|f(v)\|_{L^1([0,T],H^{s+\varepsilon}(B_+(0,1)))}$ by
$\|v\|_{L^\infty([0,T],H^{s+1}(\Omega))}$. Finally, notice that
$f(0)=0$ and thus, the Dirichlet boundary condition on $v$ naturally implies the
one on $f(v)$.
\end{demo}

\subsection{Asymptotic compactness and regularization effect}

As explained in the beginning of this section, using Duhamel formula 
$U(t)=e^{At}U_0+\int_0^t e^{A(t-s)}F(U(s))ds$ and Corollary \ref{coro-DLZ}, we
obtain two propositions related to the asymptotic smoothing effect of the
damped wave equations.
\begin{prop}\label{prop-asympt-reg}
Let $f\in\Cc^1(\RR)$ satisfying \eqref{hyp-f}, let $(u^n_0,u^n_1)$ be a
sequence of initial data which is bounded in $X=H^1_0(\Omega)\times L^2(\Omega)$
and let $(u_n)$ be the corresponding solutions of the damped wave equation
\eqref{eq}. Let $(t_n)\in\RR$ be a sequence of times such that $t_n\rightarrow
+\infty$ when $n$ goes to $+\infty$. 

Then, there exist subsequences
$(u_{\phi(n)})$ and $(t_{\phi(n)})$ and a global solution $u_\infty$ of
\eqref{eq} such that
$$\forall T>0~,~~(u_{\phi(n)},\partial_t
u_{\phi(n)})(t_{\phi(n)}+.)\longrightarrow  (u_\infty,\partial_t
u_\infty)(.)~~\text{ in } \Cc^0([-T,T],X)~.$$
\end{prop}
\begin{demo}
We use the notations of Section \ref{sect-basic}.
Due to the equivalence between the norm of $X$ and the energy given by
\eqref{equiv-norm} and the fact that the energy is decreasing in time, we know
that $U_n(t)$ is uniformly bounded in $X$ with respect to $n$
and $t\geq 0$. So, up to taking a subsequence, it weakly converges to a limit
$U_{\infty}(0)$ which gives a global solution $U_{\infty}$. We notice that, 
due to the continuity of the Cauchy problem with
respect to the initial data stated in Theorem \ref{th-Cauchy}, it is sufficient
to show that 
$U_{\phi(n)}(t_{\phi(n)})\rightarrow U_\infty(0)$ for some subsequence
$\phi(n)$. We have
\begin{align}
U_n(t_n)&= e^{At_n}U_n(0)+ \int_{0}^{t_n} e^{sA}F(U_n(t_n-s))~ds\nonumber\\
&=e^{At_n}U_n(0)+\sum_{k=0}^{\lfloor t_n\rfloor-1}e^{kA} \int_0^1 e^{sA}
F(U_n(t_n-k-s))~ds+\int_{\lfloor t_n\rfloor}^{t_n} e^{sA}F(U_n(t_n-s))~ds
\nonumber\\
&= e^{At_n}U_n(0)+ \sum_{k=0}^{\lfloor t_n\rfloor-1}
e^{kA} I_{k,n} + I_n \label{eq-asympt-reg-2}
\end{align}
Theorem \ref{th-Cauchy} shows that the Strichartz
norms $\|u_n(t_n-k-.)\|_{L^q([0,1],L^r(\Omega))}$ are uniformly bounded since
the energy of $U_n$ is uniformly bounded. Therefore, Corollary \ref{coro-DLZ}
and Proposition \ref{prop-dec-expo} show that the terms $I_{n,k}=\int_0^1 e^{sA}
F(U_n(t_n-k-s))~ds$, as well as $I_n$, are bounded by some constant $M$
in $H^{1+\varepsilon}(\Omega)\times
H^\varepsilon(\Omega)$ uniformly in $n$ and $k$. Using Proposition
\ref{prop-dec-expo} again and summing up, we get that the last terms of
\eqref{eq-asympt-reg-2} are bounded in $H^{1+\varepsilon}(\Omega)\times
H^\varepsilon(\Omega)$ uniformly in $n$ by
$$ \left\|\sum_{k=0}^{\lfloor t_n\rfloor-1}
e^{kA} I_{k,n} + I_n\right\|_{X^\varepsilon}\leq \sum_{k=0}^{\lfloor
t_n\rfloor-1} Ce^{-\lambda k}M  + M \leq M\left(1+\frac
C{1-e^{-\lambda}}\right)~. $$
Moreover, Proposition \ref{prop-dec-expo} shows that $e^{At_n}U_n(0)$ goes to
zero in $X$ when $n$ goes to $+\infty$. Therefore, by a diagonal extraction
argument and Rellich Theorem, we can extract a subsequence
$U_{\phi(n)}(t_{\phi(n)})$ that converges
to $U_\infty(0)$ in $H^1_0(B)\times L^2(B)$ for all bounded set $B$ of
$\Omega$.

To finish the proof of Proposition \ref{prop-asympt-reg}, we have to show that
this convergence holds in fact in $X$ and not only locally. Let $\eta>0$ be
given. Let $T>0$ and let 
$\tilde U_n$ be the solution of \eqref{eq} with $\tilde U_n(0)=U_n(t_n-T)$
and with $\gamma$ being replaced by $\tilde\gamma$, where
$\tilde\gamma(x)\equiv\gamma(x)$ for large $x$ and $\tilde\gamma\geq\alpha>0$
everywhere. By Proposition \ref{prop-positive-damping}, $\|\tilde
U_n(T)\|_X\leq \eta$ if $T$ is chosen sufficiently large and if $n$ is large
enough so that $t_n-T>0$. Since the information propagates at finite speed in
the wave equation, $U_n(t_n)\equiv\tilde U_n(T)$ outside a large enough bounded set and thus
$U_{\phi(n)}(t_{\phi(n)})$ has a $X-$norm smaller than $\eta$ outside this
bounded set. On the other hand, we can assume that the norm of $U_\infty(0)$ is
also smaller than $\eta$ outside the bounded set. Then, choosing $n$ large
enough, $\|U_{\phi(n)}(t_{\phi(n)})-U_\infty(0)\|_X$ becomes smaller than
$3\eta$.
\end{demo}

The trajectories $U^\infty$ appearing in Proposition
\ref{prop-asympt-reg} are trajectories which are bounded in $X$ for all times
$t\in\RR$. The following result shows that these special trajectories are more regular
than the usual trajectories of the damped wave equation. 
\begin{prop}\label{prop-reg-1}
Let $f\in\Cc^1(\RR)$ satisfying \eqref{hyp-f} and let $E_0\geq 0$. There exists
a constant $M$ such that if $u$
is a solution of \eqref{eq}, which exists for all times $t\in\RR$ and satisfies
$\sup_{t\in\RR} E(u(t))\leq E_0$, then $t\mapsto
U(t)=(u(t),\partial_tu(t))$ is continuous from $\RR$ into $D(A)$ and 
$$\sup_{t\in\RR} \|(u(t),\partial_tu(t))\|_{D(A)}\leq M~.$$
In addition, $M$ depends only on $E_0$ and the constants in
\eqref{hyp-f}.  
\end{prop}
\begin{demo}
We use a bootstrap argument. For any $t\in\RR$ and $n\in\NN$, 
$$U(t)=e^{nA}U(t-n)+\sum_{k=0}^{n-1}e^{kA} \int_0^1 e^{sA}
F(U(t-k-s))~ds~.$$
Using Proposition \ref{prop-dec-expo}, when $n$ goes to $+\infty$, we get
\begin{equation}\label{eq-reg-1}
U(t)=\sum_{k=0}^{+\infty}e^{kA} \int_0^1 e^{sA}
F(U(t-k-s))~ds~.
\end{equation}
Moreover, arguing exactly as in the proof of Proposition
\ref{prop-asympt-reg}, we show that Proposition \ref{prop-dec-expo} and
Corollary \ref{coro-DLZ} imply that Equality \eqref{eq-reg-1} also holds in
$X^\varepsilon$. Hence, $U(t)$ is uniformly bounded in $X^\varepsilon$. Then,
using again Proposition \ref{prop-dec-expo} and Corollary \ref{coro-DLZ}, 
\eqref{eq-reg-1} also holds in $X^{2\varepsilon}$ etc. Repeating the arguments
and noting that, until the last step, $\varepsilon$ only depends on $p$,
we obtain that $U(t)$ is uniformly bounded in $X^1=D(A)$.

Since the constant $C$ of Corollary \ref{coro-DLZ} only depends on
$f$ through Estimate \eqref{hyp-f}, the same holds for the bound $M$ here.
\end{demo}

\begin{prop}\label{prop-borne}
The Sobolev embedding $H^2(\Omega)\hookrightarrow \Cc^0(\Omega)$ holds and
there exists a constant $\Kc$ such that 
$$\forall u\in H^2(\Omega)~,~~\sup_{x\in\Omega}|u(x)|\leq \Kc \|u\|_{H^2}~.$$
In particular, the solution $u$ in the statement of Proposition
\ref{prop-reg-1} belongs to $\Cc^0(\overline\Omega\times\RR,\RR)$ and
$\sup_{(x,t)\in\overline\Omega\times\RR} |u(x,t)|\leq \Kc M$.
\end{prop}
\begin{demo}
Proposition \ref{prop-borne} follows directly from the fact that $\Omega$ has a
bounded geometry and from the classical Sobolev embedding
$H^2\hookrightarrow \Cc^0$ in the ball $B(0,1)$ of $\RR^3$.  
\end{demo}

%%%%%%%%%%%%%%%%%%%%%%%%%%%%%%%%%%%%%%%%%%%%%%%%%%%%%%%%%%%%%%%%%%%%%%%%%%%%
%%%%%%%%%%%%%%%%%%%%%%%%%%%%%%%%%%%%%%%%%%%%%%%%%%%%%%%%%%%%%%%%%%%%%%%%%%%%

\bigskip      

\section{Smoothness and uniqueness of non-dissipative complete
solutions}\label{sect-HR}

In this section, we consider only a non-dissipative complete solution, that is a
solution $u^*$ existing for all times $t\in\RR$ for which the energy $E$ is
constant. In other words, $u^*(t)$ solves
\begin{equation}\label{eq*}
\left\{\begin{array}{ll}
\partial^2_{tt} u^* =\Delta
u^*-\beta u^*-f(u^*)~~&(x,t)\in\Omega\times\RR~,\\
u^*(x,t)=0&(x,t)\in\partial\Omega\times \RR\\
\partial_t u^*(x,t)=0&(x,t)\in \supp(\gamma)\times \RR\\
\end{array}\right.
\end{equation}
Since the energy $E$ is not dissipated by $u^*(t)$, we can write $E(u^*)$
instead of $E(u^*(t))$. Yet, an interesting fact that will be used several times
in the sequel is that such $u^*$ is, at the same time, solution of both damped
and undamped equations.  

The purpose of this section is:\\
- First show that $u^*$ is analytic in time and smooth in space. The central
argument is to use a theorem of J.K. Hale and G. Raugel in \cite{Hale-Raugel}.\\
- Then use the unique continuation result of L. Robbiano and C. Zuily stated in
Corollary \ref{coro-RZ} to show that $u^*$ is necessarily an equilibrium point
of \eqref{eq}.\\
- Finally show that the assumption $sf(s)\geq 0$ imply
that $u^*\equiv 0$.

We enhance that the first two steps are valid and very helpful in a more general
framework than the one of our paper.

\subsection{Smoothness and partial analyticity of $u^*$}

First, we recall here the result of Section 2.2 of \cite{Hale-Raugel}, adapting
the statement to suit our notations. 
\begin{theorem}\label{th-HR}
{\bf Hale and Raugel (2003)}\\
Let $Y$ be a Banach space. Let $P_n\in\Lc(Y)$ be a sequence of continuous
linear maps and let $Q_n=Id-P_n$. Let $A:D(A)\rightarrow Y$ be the generator of
a continuous semigroup $e^{tA}$ and let $G\in\Cc^1(Y)$. We assume that $V$ is a
complete mild solution in $Y$ of  
$$\partial_t V(t)=A V(t)+G(V(t))~~~\forall t\in\RR~.$$
We further assume that\\[2mm]
(i) $\{V(t),t\in\RR\}$ is contained in a compact set $K$ of $Y$.\\
(ii) for any $y\in Y$, $P_n y$ converges to $y$ when
$n$ goes $+\infty$ and $(P_n)$ and $(Q_n)$ are sequences of $\Lc(Y)$ bounded
by $K_0$.\\
(iii) the operator $A$ splits in $A=A_1+B_1$ where $B_1$ is bounded and $A_1$
commutes with $P_n$.\\ 
(iv) there exist $M$ and $\lambda>0$ such that $\|e^{At}\|_{\Lc(Y)}\leq
M e^{-\lambda t}$ and $\|e^{(A_1+Q_nB_1)t}\|_{\Lc(Q_nY,Y)}\leq
M e^{-\lambda t} $ for all $t\geq 0$.\\
(v) $G$ is analytic in the ball $B_Y(0,r)$, where $r$ is such that $r\geq 4K_0
\sup_{t\in\RR} \|V(t)\|_Y$. More precisely, there exists $\rho>0$ such that $G$
can be extended to an holomorphic function of $B_Y(0,r)+i B_Y(0,\rho)$.\\
(vi) $\left\{DG(V(t))V_2\left|t\in \R, \nor{V_2}{Y}\leq 1 \right.\right\}$ is a
relatively compact set of $Y$.\\[2mm]
Then, the solution $V(t)$ is analytic from $t\in\RR$ into $Y$.
\end{theorem}

More precisely, Theorem \ref{th-HR} is Theorem 2.20 (which relates to Theorem
2.12) of \cite{Hale-Raugel} applied with Hypothesis (H3mod) and (H5). 
 
Proposition \ref{prop-reg-1} shows that $u^*$ is continuous in both space and
time variables. We apply Theorem \ref{th-HR} to show that, because $f$ is
analytic, $u^*$ is also analytic with respect to the time. 
\begin{prop}\label{prop-reg-time}
Let $f\in\Cc^1(\RR)$ satisfying \eqref{hyp-f} and let $E_0\geq 0$.
Let $\Kc$ and $M$ be the constants given by Propositions \ref{prop-reg-1} and
\ref{prop-borne}. Assume that $f$
is analytic in $[-4\Kc M,4\Kc M]$. Then for any non-dissipative complete
solution $u^*(t)$ solving \eqref{eq*} and satisfying $E(u^*)\leq E_0$,
$t\mapsto u^*(.,t)$ is analytic from $\RR$ into $X^\alpha$
with $\alpha\in (1/2,1)$. In particular, for
all $x\in\Omega$, $u^*(x,t)$ is analytic with respect to the time.
\end{prop}
\begin{demo}
Theorem \ref{th-HR} uses strongly some compactness properties. Therefore,
we need to truncate our solution to apply the theorem on a bounded domain (of
course, this is not necessary and easier if $\Omega$ is already bounded). 

Let 
$\chi\in\Cc^\infty_0(\overline \Omega)$ be such that
$\Drond{\chi}{\nu}=0$ on $\partial\Omega$, $\chi\equiv 1$ in
$\{x\in\Omega,\gamma(x)=0\}$ and $\supp(\chi)$ is included in a smooth bounded
subdomain $\Oc$ of $\Omega$. Since Proposition \ref{prop-reg-1} shows that
$u^*\in\Cc^0(\RR,D(A))$ and since $u^*$ is constant with respect to the time in
$\supp(\gamma)$, $(1-\chi)u^*$ is obviously analytic from $\RR$ into $D(A)$. It
remains to obtain the analyticity of $\chi u^*$.

In this proof, the damping $\gamma$ needs to be more regular than just
$L^\infty(\Omega)$. We replace $\gamma$ by a damping $\tilde\gamma\in
C^{\infty}(\Omega)$, which has the same geometrical properties (GCC) and
\eqref{hyp-gamma} and which vanishes where $\gamma$ does. Notice that $\gamma
\partial_t u^*\equiv 0 \equiv \tilde \gamma\partial_t u^*$, therefore
replacing $\gamma$ by $\tilde \gamma$ has no consequences here.

Let $v=\chi u^*$, we have
\begin{equation}\label{eq-v*}
\left\{\begin{array}{ll}
\partial^2_{tt} v + \tilde\gamma(x) \partial_t v =\Delta
v-\beta v +g(x,v)~~&(x,t)\in \Oc\times\RR_+~,\\
v(x,t)=0&(x,t)\in\partial\Oc\times \RR_+\\
\end{array}\right.
\end{equation}
with $g(x,v)=-\chi(x)f(v+(1-\chi)u^*(x))-2(\grad\chi.\grad
u^*)(x)-(u^*\Delta\chi)(x)$. We apply Theorem \ref{th-HR} with the following
setting. Let
$Y=X^\alpha=H^{1+\alpha}(\Oc)\cap H^1_0(\Oc)\times H^\alpha_0(\Oc)$
with $\alpha\in (1/2,1)$. Let $V=(v,\partial_t v)$ and let $G(v)=(0,g(.,v))$. We
set 
$$A= A_1+B_1=\left(\begin{array}{cc}
0& Id\\ \Delta-\beta&0 \end{array}\right)+ \left(\begin{array}{cc}
0& 0\\ 0&-\tilde \gamma \end{array}\right)~.$$
Let $(\lambda_k)_{k\geq 1}$ be the negative eigenvalues of
the Laplacian operator on $\Oc$ with Dirichlet boundary conditions and let
$(\varphi_k)$ be corresponding eigenfunctions. We set $P_n$ to be the canonical
projections of $X$ on the subspace generated by $((\varphi_k,0))_{k=1...n}$ and
$((0, \varphi_k))_{k=1...n}$.

To finish the proof of Proposition \ref{prop-reg-time}, we only have to check
that the hypotheses of Theorem \ref{th-HR} hold. 

The trajectory $V$ is compact since we know by Proposition \ref{prop-reg-1}
 that it is bounded in $X^1$, which gives (i).
 
Hypothesis (ii) and (iii) hold with $K_0=1$
by construction of $P_n$ and because $B_1$ is bounded in $Y$ since $\tilde
\gamma$ belongs to $C^{\infty}(\Omega)$. 

The first part of Hypothesis (iv) follows from Proposition
\ref{prop-dec-expo}. The second estimate
$\|e^{(A_1+Q_nB_1)t}\|_{\Lc(Q_nY,Y)}\leq
M e^{-\lambda t}$ means that the restriction $Q_nAQ_n$ of $A$ to the
high frequencies of the wave operator also generates a semigroup satisfying an
exponential decay. By a result of Haraux in \cite{Haraux3}, we know that this
also holds, see Section 2.3.2 of \cite{Hale-Raugel} for the detailed arguments.
 
We recall that $u^*(x,.)$ is constant outside
$\chi^{-1}(1)$ and belongs locally to $H^{1+\alpha}$ since $u^*\in D(A)$.
Therefore, the terms $(1-\chi)u^*(x)$, $\grad\chi.\grad
u^*$ and $u^*\Delta\chi$ appearing in the definition of $g$ are in $H^{1}$.
Moreover, they satisfy Dirichlet boundary condition on $\partial\Omega$ since
$u^*\equiv 0$ and $\partial_\nu \chi\equiv 0$ there. Of course, they
also satisfy Dirichlet boundary condition on the other parts of $\partial\Oc$
since $\chi\equiv 0$ outside $\Oc$. Notice that
$\alpha>1/2$ and thus $H^{1+\alpha}(\Oc)\cap H^1_0(\Oc)$ is an algebra
included in $\Cc^0$. Therefore, \eqref{hyp-f} shows that $G$ is of
class $\Cc^1$ in the bounded sets of $Y$. Since 
$u\in[-4\Kc M,4\Kc M] \mapsto f(u)\in\RR$ is analytic, it can be extended to a
holomorphic function in $[-4\Kc M,4\Kc M]+i[-\rho,\rho]$ for small $\rho>0$.
Using again the embedding $H^{1+\alpha}(\Oc)\hookrightarrow \Cc^0(\Oc)$ and the
definitions of $\Kc$ and $M$, we deduce that (v) holds.

Finally, for $V_2=(v_2,\partial_t v_2)$ with $\nor{V_2}{Y}\leq 1$
$DG(V(t))V_2=(0,-\chi(x)f'(v(t)+(1-\chi)u^*(x))v_2)$ is relatively compact in
$Y$ since $v(t)$ is bounded in $H^2\cap H^1_0$ due to Proposition
\ref{prop-reg-1} and therefore $v_2\in H^{1+\alpha}\mapsto
\chi(x)f'(v(t)+(1-\chi)u^*(x))v_2 \in H^\alpha$ is a compact map. This yields
(vi).
\end{demo}

Once the time-regularity of $u^*$ is proved, the space-regularity follows
directly.
\begin{prop}\label{prop-reg-space}
Let $f$ and $u^*$ be as in Proposition \ref{prop-reg-time}, then
$u^*\in\Cc^\infty(\Omega\times\RR)$.
\end{prop}
\begin{demo}
Proposition \ref{prop-reg-time} shows that $u^*$ and all its time-derivatives
belongs to $X^\alpha$
with $\alpha\in (1/2,1)$. Due to the Sobolev embeddings, this implies that any
time-derivative of $u^*$ is H\"older continuous. Writing 
\begin{equation}\label{eq-reg-space}
\Delta u^*=\partial^2_{tt}u^*+\beta u^*+f(u^*)
\end{equation}
and using the local elliptic regularity properties (see \cite{Miranda} and the
references therein for example), we get that $u^*$ is locally of class
$\Cc^{2,\lambda}$ in space for some $\lambda\in(0,1)$. Thus, $u^*$ is of class
$\Cc^{2,\lambda}$ in both time and space. Then, we can use a bootstrap argument
in \eqref{eq-reg-space} to show that $u^*$ is of class $\Cc^{2k,\lambda}$ for
all $k\in\NN$.
\end{demo}

\subsection{Identification of $u^*$}

The smoothness and the partial analyticity of $u^*$ shown in Propositions
\ref{prop-reg-time} and \ref{prop-reg-space} enable us to use the unique
continuation result of \cite{Rob-Zui}.
\begin{prop}\label{prop-eq}
Let $f$ and $u^*$ be as in Proposition \ref{prop-reg-time}, then $u^*$ is
constant in time, i.e. $u^*$ is an equilibrium point of
the damped wave equation \eqref{eq}.
\end{prop}
\begin{demo}
Setting $v=\partial_t u^*$, we get $\partial^2_{tt}v=\Delta v-\beta v - f'(u^*)
v$. Propositions \ref{prop-reg-time} and \ref{prop-reg-space} show that $u^*$ is
smooth and analytic with respect to the time and moreover $v\equiv 0$ in
$\supp(\gamma)$. Thus, the unique continuation result stated in Corollary
\ref{coro-RZ} yields $v\equiv 0$ everywhere. 
\end{demo}

The sign assumption on $f$ directly implies that $0$ is the only possible
equilibrium point of \eqref{eq}.
\begin{corollary}\label{coro-uniNL}
Let $f\in\Cc^1(\RR)$ satisfying \eqref{hyp-f} and let $E_0\geq 0$.
Let $\Kc$ and $M$ be the constants given by Propositions \ref{prop-reg-1} and
\ref{prop-borne} and assume that $f$ is analytic in $[-4\Kc M,4\Kc M]$. Then,
the unique solution $u^*$ of (\ref{eq*}) with $E(u^*)\leq E_0$ is $u^*\equiv 0$.
\end{corollary}
\begin{demo}
Due to Proposition \ref{prop-eq}, $u^*$ is solution of $\Delta
u^*-\beta u^*=f(u^*)$. By multiplying by $u^*$ and integrating by part, we
obtain $\int_{\Omega} |\nabla u^*|^2+\beta|u^*|^2\, dx=-\int_{\Omega}
u^*f(u^*)\, dx$, which is non-positive due to Assumption \eqref{hyp-f}. Since
$\beta\geq 0$ is such that $\Delta-\beta$ is negative
definite, this shows that $u^*\equiv 0$.
\end{demo}

%%%%%%%%%%%%%%%%%%%%%%%%%%%%%%%%%%%%%%%%%%%%%%%%%%%%%%%%%%%%%%%%%%%%%%%%%%%%
%%%%%%%%%%%%%%%%%%%%%%%%%%%%%%%%%%%%%%%%%%%%%%%%%%%%%%%%%%%%%%%%%%%%%%%%%%%%

\bigskip      

\section{Proof of Theorem \ref{th1}}\label{sect-th1}

Due to Proposition \ref{critere-decay}, Theorem \ref{th1} directly
follows from the following result.
\begin{prop}\label{thminegobserv}
Let $f\in\Cc^1(\RR)$ satisfying \eqref{hyp-f} and let $E_0\geq 0$.
Let $\Kc$ and $M$ be the constants given by Propositions \ref{prop-reg-1} and
\ref{prop-borne}. Assume that $f$ is analytic in $[-4\Kc M,4\Kc M]$ and that
$\gamma$ is as in Theorem \ref{th1}. Then, there exist $T>0$ and $C>0$ such that
for any $u$ solution of \eqref{eq} with $E(u)(0)\leq E_0$ satisfies 
$$
E(u)(0) \leq C\iint_{[0,T]\times \Omega}  \gamma(x)\left|\partial_t u\right|^2 ~dtdx .
$$
\end{prop}
\begin{demo} 
We argue by contradiction: we assume that there exists a sequence $(u_n)$ of
solutions of \eqref{eq} and a sequence of times $(T_n)$ converging to
$+\infty$ such that
\begin{equation}\label{eq-demo-1}
\iint_{[0,T_n]\times \Omega}  \gamma(x)\left|\partial_t u_n\right|^2 ~dtdx 
\leq \frac{1}{n}E(u_n)(0)\leq \frac{1}{n}E_0.
\end{equation}
Denote $\alpha_n=(E(u_n)(0))^{1/2}$. Since $\alpha\in[0,\sqrt{E_0}]$, we can
assume that $\alpha_n$ converges to a limit $\alpha$ when $n$ goes to
$+\infty$. We distinguish two cases: $\alpha>0$ and $\alpha=0$.

\bigskip

\textbullet \quad First case: $\alpha_n \longrightarrow \alpha>0$\\
Notice that, due to \eqref{equiv-norm}, $\|(u_n,\partial_t u_n)(0)\|_X$ is
uniformly bounded from above and from below by positive numbers. We set
$u^*_n=u_n(T_n/2+.)$. Due to the asymptotic compactness property stated in
Proposition \ref{prop-asympt-reg}, we can assume that $u^*_n$ converges to a
solution $u^*$ of \eqref{eq} in $\Cc^0([-T,T],X)$ for all time $T>0$. We notice
that 
$$E(u_n(0))\geq E(u^*_n(0))=E(u_n(0))-\iint_{[0,T_n/2]\times \Omega} 
\gamma(x)\left|\partial_t u_n\right|^2 \geq (1-1/n) E(u_n(0))$$
and thus $E(u^*(0))=\alpha^2>0$. Moreover, \eqref{eq-demo-1} shows that
$\gamma\partial_t u^*_n$ converges to zero in $L^2([-T,T],L^2(\Omega))$ for
any $T>0$ and thus $\partial_t u^*\equiv 0$ in $\supp(\gamma)$. In other words,
$u^*$ is a non-dissipative solution of \eqref{eq}, i.e. a solution of
\eqref{eq*} with $E(u^*)=\alpha^2\leq E_0$. Corollary \ref{coro-uniNL} shows
that $u^*\equiv 0$, which
contradicts the positivity of $E(u^*(0))$.

\bigskip

\textbullet \quad Second case: $\alpha_n \longrightarrow 0$\\
The assumptions on $f$ allow to write $f(s)=f'(0)s+R(s)$ with 
\begin{equation}\label{eq-demo-3}
|R(s)|\leq C(|s|^2+|s|^p)~~\text{ and }~~|R'(s)|\leq C(|s|+|s|^{p-1})~.
\end{equation}
Let us make the change of unknown $w_n=u_n/\alpha_n$. Then, $w_n$ solves
\begin{equation}
\label{eqnalphan}
\Box w_n +\gamma(x)\partial_t w_n+(\beta+f'(0))w_n+\frac{1}{\alpha_n}R(\alpha_nw_n)=0
\end{equation}
and
\begin{equation}\label{eq-demo-2}
\iint_{[0,T_n]\times \Omega}  \gamma(x)\left|\partial_t w_n\right|^2 ~dtdx \leq
\frac{1}{n}. 
\end{equation}
Denote $W_n=(w_n,\partial_t w_n)$. Due to the equivalence between norm and
energy given by \eqref{equiv-norm}, the 
scaling $w_n=u_n/\alpha_n$ implies that
$\left\|(w_n(0),\partial_tw_n(0))\right\|_{X}$ is uniformly bounded from 
above and from below by positive numbers. Moreover, \eqref{eq-demo-1} implies 
\begin{equation}
\label{boundedbelow} 
\left\|W_n(t)\right\|_{X}=\frac{\left\|(U_n(t))\right\|_{X}}{\alpha_n}\geq
C\frac{E(u_n(t))^{1/2}}{\alpha_n}\geq C\frac{(E(u_n
)(0)-\alpha_n^2/n)^{1/2}}{\alpha_n}\geq \frac C2 >0
\end{equation}
for any $t\in [0,T_n]$ and $n$ large enough.

We set $f_n=1/\alpha_n R(u_n)$ and $F_n=(0,f_n)$. The stability estimate of
Proposition \ref{th-Cauchy} implies that $\nor{u_n}{L^q([k,k+1],L^r)}\leq C
\alpha_n $ uniformly for $n, k\in \NN$. In particular, combined with
\eqref{eq-demo-3}, this gives 
$$
\nor{f_n}{L^1([k,k+1],L^2)}=\nor{\frac{1}{\alpha_n}R(\alpha_nw_n)}{L^1([k,k+1],
L^2)}\leq C (\alpha_n+\alpha_n^{p-1}) 
$$
We can argue as in Proposition \ref{prop-asympt-reg} and write
\begin{align}
W_n(T_n)&=e^{\tilde{A}T_n}W_n(0)+\sum_{k=0}^{\lfloor
T_n\rfloor-1}e^{\tilde{A}(T_n-k)} \int_0^1 e^{-\tilde{A}s}
F_n(k+s)~ds \nonumber\\
&~~~~~~~~ +e^{ \tilde{A}(T_n-\lfloor T_n\rfloor)}\int_0^{T_n-\lfloor
T_n\rfloor} e^{-\tilde{A}s} F_n(\lfloor T_n\rfloor+s)~ds~. \label{eq-asympt-reg}
\end{align}
  where $\tilde{A}$ is the modified damped wave operator $$\tilde{A}=\left(\begin{array}{cc}
0& Id\\ \Delta-\beta-f'(0) &-\gamma \end{array}\right)~. $$ 
 Notice that $\tilde{A}$ satisfies an exponential decay as $A$ in Proposition
\ref{prop-dec-expo} since \eqref{hyp-f} implies $f'(0)\geq 0$. By summing up as
in Proposition \ref{prop-asympt-reg}, we get 
$$
 \nor{W_n(T_n)}{X}\leq Ce^{-\lambda T_n}+C (\alpha_n+\alpha_n^{p-1})
$$
 which goes to zero, in a contradiction with \eqref{boundedbelow}.
\enp

As a direct consequence of Proposition \ref{thminegobserv}, we obtain a unique
continuation property for nonlinear wave equations. Notice that the time of
observation $T$ required for the unique continuation is not explicit. Thus, this result is not so convenient as a unique continuation property. But it may be useful for other nonlinear stabilization problems as $\Box u+\gamma(x)g(\partial_t u)+f(u)=0$.
\begin{corollary} \label{corrUCP}
Let $f\in\Cc^1(\RR)$ satisfying \eqref{hyp-f} and let $E_0\geq 0$.
Assume that $f$ is analytic in $\R$ and that $\omega$ is an open
subset of $\Omega$ satisfying (GCC). Then, there exist $T>0$ such that
the only solution $u$ of 
\bneqn\label{eq-iiibis}
\Box u+\beta u+f(u)&=&0 \quad \textnormal{ on }[-T,T]\times \Omega~\\
\partial_tu&\equiv&0\quad \textnormal{ on }[-T,T]\times \omega~. 
\eneqn
with $E(u)(0)\leq E_0$ is $u\equiv 0$.
\end{corollary}
\begin{demo}
Corollary \ref{corrUCP} is a straightforward consequence of Proposition
\ref{thminegobserv} since we can easily construct a smooth damping $\gamma$
supported in $\omega$ and such that $Supp(\gamma)$ satisfies (GCC). We only have
to remark that $u$ solution of \eqref{eq-iiibis} is also solution of 
\eqref{eq}. 
\end{demo}
%%%%%%%%%%%%%%%%%%%%%%%%%%%%%%%%%%%%%%%%%%%%%%%%%%%%%%%%%%%%%%%%%%%%%%%%%%%%
%%%%%%%%%%%%%%%%%%%%%%%%%%%%%%%%%%%%%%%%%%%%%%%%%%%%%%%%%%%%%%%%%%%%%%%%%%%%

\bigskip      

\section{Proof of Theorem \ref{th2}}\label{sect-th2}
Before starting the proof of Theorem \ref{th2} itself, we prove that
$\Cg^1(\RR)$ is a Baire space, that is that any countable intersection of open
dense sets is dense. This legitimizes the
genericity in $\Cg^1(\RR)$ as a good notion of large subsets of $\Cg^1(\RR)$.
We recall that $\Cg^1(\RR)$ is defined by \eqref{eq-Cg}
and endowed by Whitney topology, the open sets of which are generated by the
neighbourhoods $\Nc_{f,\delta}$ defined by \eqref{eq-Whitney}. 
\begin{prop}\label{prop-Baire}
The space $\Cg^1(\RR)$ endowed with Whitney topology is a Baire space.
\end{prop}
\begin{demo}
The set $\Cg^1(\RR)$ is not an open set of $\Cc^1(\RR)$, and neither a
submanifold. It is a closed subset of $\Cc^1(\RR)$, but $\Cc^1(\RR)$ endowed
with Whitney topology is not a completely metrizable space, since it is not
even metrizable (the neighbourhoods of a function $f$ are not generated by a
countable subset of them). Therefore, we have to come back to the basic proof
of Baire property as in \cite{GG}.

Let $\Uc$ be an open set of $\Cg^1(\RR)$ and let $(\Oc_n)_{n\in\NN}$ be a
sequence of open dense sets of $\Cg^1(\RR)$. By density, there exists a function
$f_0\in\Cg^1(\RR)$ in  $\Uc\cap\Oc_0$ and by openness, there exists a
positive continuous function $\delta_0$ such that the neighbourhood
$\Nc_{f_0,\delta_0}$ is contained in $\Uc\cap\Oc_0$. By
choosing $\delta_0$ small enough, one can also assume that
$\Nc_{f_0,2\delta_0}\subset \Uc\cap\Oc_0$ and that
$\sup_{u\in\RR}|\delta_0(u)|\leq 1/2^0$. By recursion,
one constructs similar balls $\Nc_{f_n,\delta_n}\subset
\Nc_{f_{n-1},\delta_{n-1}}\subset \Nc_{f_{0},\delta_0}$ such that
$\Nc_{f_n,2\delta_n}\subset \Uc\cap\Oc_n$ and that
$\sup_{u\in\RR}|\delta_n(u)|\leq 1/2^n$. Since $\Cc^1([-m,m],\RR)$ endowed with
the uniform convergence topology is a complete metric space, the sequence
$(f_i)$ converges to a function $f\in\Cc^1(\RR,\RR)$ uniformly in any compact
set of $\RR$. By construction, the limit $f$ satisfies 
\begin{equation}\label{eq-demo-Baire2}
\forall n\in\NN,~\forall u\in\RR~,~~\max (|f(u)-f_n(u)|,|f'(u)-f'_n(u)|)
\leq\delta_n(u)<2\delta_n(u)~,
\end{equation}
as well as $f(0)=0$ and $uf(u)\geq 0$ since any $f_n$ satisfies \eqref{hyp-f}.
Moreover, there exist $C>0$ and $p\in[1,5)$ such that $f_0$ satisfies
\begin{equation}\label{eq-demo-Baire}
|f_0(u)|\leq C(1+|u|)^p~\text{ and }~~|f'_0(u)|\leq C(1+|u|)^{p-1}~.
\end{equation}
Since $\max (|f(u)-f_0(u)|,|f'(u)-f'_0(u)|) \leq\delta_0(u)\leq 1$, $f$ also
satisfies \eqref{eq-demo-Baire} with a constant $C'=C+1$. Therefore, $f$
satisfies \eqref{hyp-f} and thus belongs to $\Cg^1(\RR)$. In addition, $f$
satisfying \eqref{eq-demo-Baire2} and $\Nc_{f_n,2\delta_n}$ being
contained in $\Uc\cap\Oc_n$, we get $f\in \Uc\cap\Oc_n$ for all $n$. 
This shows that $\cap_{n\in\NN}\Oc_n$ intersects any open set $\Uc$ and
therefore is dense in $\Cg^1(\RR)$.
\end{demo}

{ \noindent \emph{\textbf{Proof of Theorem \ref{th2}:}}}
We denote by $\Gg_n$ the set of functions $f\in\Cg^1(\RR)$ such that the
exponential decay property (ED) holds for $E_0=n$. Obviously,
$\Gg=\cap_{n\in\NN}\Gg_n$ and hence it is sufficient to prove that $\Gg_n$ is an
open dense subset of $\Cg^1(\RR)$. We sketch here the main arguments to prove
this last property.

\vspace{3mm}

\noindent{\bf $\Gg_n$ is a dense subset: }let $\Nc$ be a neighbourhood of
$f_0\in\Cg^1(\RR)$. Up to choosing $\Nc$ smaller, we can assume that the
constant in \eqref{hyp-f} is independent of $f\in\Nc$. Due to Propositions
\ref{prop-reg-1} and \ref{prop-borne}, there exist constants $\Kc$ and $M$ such
that, for all $f\in\Nc$, all the global non-dissipative trajectories $u$ of
\eqref{eq} with $E(u)\leq n$ are such that
$\|u\|_{L^\infty(\Omega\times\RR)}\leq \Kc M$. We claim that we 
can choose
$f\in\Nc$ as close to $f_0$ as wanted such that $f$ is analytic on $[-4\Kc
M,4\Kc M]$ and still satisfies \eqref{hyp-f}. Then, Proposition
\ref{thminegobserv} shows that $f$ satisfies (ED) with $E_0=n$ i.e. that
$f\in\Gg_n$.  

To obtain this suitable function $f$, we proceed as follows. First, we set
$a=4\Kc M$ and notice that it is sufficient to explain how we construct $f$ in
$[-a,a]$. Indeed, one can easily extend a perturbation $f$ of $f_0$ in $[-a,a]$
satisfying $f(s)s\geq 0$ to a perturbation $\tilde f$ of $f_0$ in $\RR$, equal
to $f_0$ outside of $[-a-1,a+1]$ and such that $f(s)s\geq 0$ in $[-a-1,a+1]$. We
construct $f$ in $[-a,a]$ as follows. Since $f_0(s)s\geq 0$, we have that
$f'_0(0)\geq 0$. We perturb $f_0$ to $f_1$ such that $f_1(0)=0$, $f'_1(s)\geq
\varepsilon>0$ in a small interval $[-\eta,\eta]$ and $sf_1(s)\geq 2\varepsilon$
in $[-a,-\eta]\cup [\eta,a]$, where $\varepsilon$ could be chosen as small as
needed. Then we perturb $f_1$ to obtain a function $f_2$ which is analytic in
$[-a,a]$ and satisfies $f'_2(s)>0$ in $[-\eta,\eta]$, $sf_2(s)\geq\varepsilon$
in $[-a,-\eta]\cup [\eta,a]$ and $|f_2(0)|<\varepsilon / a$. Finally, we set
$f(s)=f_2(s)-f_2(0)$ and check that $f$ is analytic and satisfies $sf(s)\geq 0$
in $[-a,a]$. Moreover, up to choosing $\varepsilon$ very small, $f$ is as close
to $f_0$ as wanted. 

\vspace{3mm}

\noindent{\bf $\Gg_n$ is an open subset: } let $f_0\in\Gg_n$. Proposition
\ref{prop-dec-expo} shows the existence of a
constant $C$ and a time $T$ such that, for all solution $u$ of \eqref{eq}, 
\begin{equation}\label{eq-demo-th2}
E(u(0))\leq E_0~~\Longrightarrow~~ E(u(0))\leq C\int_0^T\int_\Omega
\gamma(x)|\partial_t u(x,t)|^2~dxdt~.
\end{equation}
The continuity of the trajectories in $X$ with respect to
$f\in\Cg^1(\RR)$ is not difficult to obtain: using the strong control of $f$
given by Whitney topology, the arguments are the same as the
ones of the proof of the continuity with respect to the initial data, stated in
Theorem \ref{th-Cauchy}. Thus,
\eqref{eq-demo-th2} holds 
also for any $f$ in a neighbourhood $\Nc$ of $f_0$, replacing the constant $C$
by 
a larger one. Therefore, Proposition \ref{prop-dec-expo} shows that
$\Nc\subset\Gg_n$ and hence that $\Gg_n$ is open.
{\hfill$\square$\\ \vspace{0.4cm}}

%%%%%%%%%%%%%%%%%%%%%%%%%%%%%%%%%%%%%%%%%%%%%%%%%%%%%%%%%%%%%%%%%%%%%%%%%%%%
%%%%%%%%%%%%%%%%%%%%%%%%%%%%%%%%%%%%%%%%%%%%%%%%%%%%%%%%%%%%%%%%%%%%%%%%%%%%

\bigskip      \section{A proof of compactness and regularity with the usual
arguments of control theory}\label{sect-control}

In this section, we give an alternative proof of the compactness and
regularity properties of Propositions \ref{prop-asympt-reg} and \ref{prop-reg-1}.
We only give its outline
since it is redundant with previous results of the article. Moreover, it
is quite similar to the arguments of \cite{DLZ}. Yet, the arguments of this
section are interesting because they do not require any asymptotic arguments
and they show a regularization effect through an observability estimate with a
finite time $T$, which can be explicit. However, for the moment, it seems
impossible to obtain an analytic regularity similar to Proposition
\ref{prop-reg-time} with this kind of arguments. 

Instead of using a Duhamel formula with an infinite interval of time $(-\infty,t)$ as in \eqref{eq-reg-1},
the main idea is to use as a black-box an observability estimate for $T$ large
enough, $T$ being the time of geometric control condition, 
\begin{equation}
\label{inegobserv}
\nor{U_0}{X^s}^2 \leq C\nor{B e^{tA}U_0}{L^2([0,T],X^s)}^2
\end{equation}
where 
$$
A=\left(\begin{array}{cc}
0& Id\\ \Delta-\beta &-\gamma \end{array}\right)~ \textnormal{ and } B=\left(\begin{array}{cc}
0& 0\\ 0&-\gamma \end{array}\right)~. $$
 
The first aim is to prove that a solution of \eqref{eq*}, globally bounded in
energy, is also globally bounded in $X^s$ for $s\in [0,1]$. We proceed step by
step. First, let us show that it is bounded in $X^\varepsilon$.
\begin{itemize}
 \item we fix $T>$ large enough to get the observability estimate
\eqref{inegobserv}. By the existence theory on each $[t_0,t_0+T]$,
$u_{|[t_0,t_0+T]}$ is
bounded in Strichartz norms, uniformly for $t_0\in \RR$. Since the nonlinearity
is subcritical, Corollary \ref{coro-DLZ} gives that $f(u)$ is globally bounded
in $L^1([t_0,t_0+T],H^{1+\varepsilon})$. 
\item we decompose the solution into its linear and nonlinear part by the Duhamel formula 
$$
U(t)=e^{A(t-t_0)}U(t_0)+\int_{t_0}^t
e^{A(t_0-\tau)}f(U(\tau))d\tau=U_{lin}+U_{Nlin}.
$$
Since $f(u)$ is bounded in $L^1([{t_0},{t_0}+T],H^{1+\varepsilon})$, $U_{Nlin}$
is uniformly bounded in $C([{t_0},{t_0}+T],X^{\e})$.
\item we will now use the linear observability estimate \eqref{inegobserv} with $s=\e$, applying it to $U_{lin}$:
\begin{equation}
\label{inegdecomp}
\nor{U({t_0})}{X^{\e}}^2=\nor{U_{lin}({t_0})}{X^{\e}}^2 \leq C
\int_{t_0}^{{t_0}+T}
\nor{\gamma(x)\partial_t u_{lin} }{H^{\varepsilon}}^2. 
\end{equation}
Then, using triangular inequality, we get 
\bna 
\int_{t_0}^{{t_0}+T} \nor{\gamma(x)\partial_t u_{lin} }{H^{\varepsilon}}^2&\leq
&2\int_{t_0}^{{t_0}+T} \nor{\gamma(x)\partial_t u
}{H^{\varepsilon}}^2+2\int_{t_0}^{{t_0}+T} \nor{\gamma(x)\partial_t u_{Nlin}
}{H^{\varepsilon}}^2\\
&\leq &2\int_{t_0}^{{t_0}+T} \nor{\gamma(x)\partial_t u_{Nlin}
}{H^{\varepsilon}}^2 \leq C
\ena
where we have used $\partial_t u\equiv 0$ on $\omega$ and $U_{Nlin}$ is bounded
in $C([{t_0},{t_0}+T],X^{\e})$. Combining this with \eqref{inegdecomp} for any
$t_0\in\RR$, we obtain that $U$ is uniformly bounded in $X^{\e}$ on $\R$.
\end{itemize}
Repeating the arguments, we show that $u$ is bounded in $X^{2\varepsilon}$,
$X^{3\varepsilon}$ and so on\ldots until $X^1$. Similar ideas allows to prove
some theorem of propagation of compactness in finite time, replacing the
asymptotic compactness property of Proposition \ref{prop-asympt-reg}.

As said above, an advantage of this method, compared to the one used in
Propositions \ref{prop-asympt-reg} and \ref{prop-reg-1}, is that it allows to propagate the regularity or
the compactness on some finite interval of fixed length. Yet, it seems that such
propagation results are not available in the analytic setting. Indeed, it seems
that, for nonlinear equations, the propagation of analytic  regularity or of
nullity in finite time is much harder to prove. We can for
instance refer to the weaker (with respect to the geometry) result of
Alinhac-M\'etivier \cite{Alin_Met} or the negative result of M\'etivier
\cite{Met_contrex}.
%%%%%%%%%%%%%%%%%%%%%%%%%%%%%%%%%%%%%%%%%%%%%%%%%%%%%%%%%%%%%%%%%%%%%%%%%%%%
%%%%%%%%%%%%%%%%%%%%%%%%%%%%%%%%%%%%%%%%%%%%%%%%%%%%%%%%%%%%%%%%%%%%%%%%%%%%

\bigskip      \section{Applications}\label{sect-disc}
\subsection{Control of the nonlinear wave equation}

In this subsection, we give a short proof of Theorem \ref{thmcontrol}, which states the global
controllability of the nonlinear wave equation. The first step consists in a local control theorem. 
\begin{theorem}[Local control] 
\label{thmlocalcontrol}
Let $\omega$ satisfying the geometric control condition for a time $T$. Then, there
exists $\delta$ such that for any $(u_0,u_1)$ in $H^1_0(\Omega)\times
L^2(\Omega)$, with  
$$
\nor{(u_0,u_1)}{H^1_0\times L^2} \leq \delta
$$
there exists $g\in L^{\infty}([0,T],L^2)$ supported in $[0,T]\times \omega$ such
that the unique strong solution of  
\begin{eqnarray*}
\left\lbrace
\begin{array}{rcl}
\Box u+\beta u+f(u)&=& g\quad \textnormal{on}\quad [0,T]\times \Omega\\
(u(0),\partial_t u(0))&=&(u_0,u_1) .
\end{array}
\right.
\end{eqnarray*}
satisfies $(u(T),\partial_t u(T))=(0,0)$.
\end{theorem}
\begin{demo}
The proof is exactly the same as Theorem 3 of \cite{DLZ} or Theorem 3.2 of
\cite{LaurentNLW}. The main argument consists in seeing the problem as a
perturbation of the linear controllability which is known to be true in  our
setting. 
\enp

\medskip

Now, as it is very classical, we can combine the local controllability with our
stabilization theorem to get global controllability. 

\medskip

{ \noindent \emph{\textbf{Sketch of the proof of Theorem \ref{thmcontrol}:}}}
In a first step, we choose as a control $g=-\gamma(x)\partial_t \widetilde{u}$ where $\widetilde{u}$
is solution of (\ref{eq}) with initial data $(u_0,u_1)$. By uniqueness of
solutions, we have $u=\widetilde{u}$. Therefore, thanks to Theorem \ref{th1},
for a large time $T_1$, only depending on $R_0$, we have $\nor{(u(T_1),\partial_t
u(T_1))}{H^1\times L^2} \leq \delta$. Then, Theorem
\ref{thmlocalcontrol} allows to find a control that brings $(u(T_1),\partial_t
u(T_1))$ to $0$. In other words, we have found a control $g$
supported in $\omega$ that brings $(u_0,u_1)$ to $0$.
We obtain the same result for $(\tilde{u}_0,\tilde{u}_1)$ and
conclude, by reversibility of the equation, that we can also bring $0$ to $(\tilde{u}_0,\tilde{u}_1)$. 
{\hfill$\square$\\ \vspace{0.4cm}}

\subsection{Existence of a compact global attractor}
In this subsection, we give the modification of the proofs of this paper
necessary to
get Theorem \ref{th-attrac} about the existence of a global attractor. 

 The energy
associated to \eqref{eq-attrac} in $X=H^1_0(\Omega)\times L^2(\Omega)$ is given
by 
$$E(u,v)=\int_{\Omega} \frac 12(|\grad u|^2+|v|^2)+V(x,u)dx~,$$
where $V(x,u)=\int_0^u f(x,\xi)d\xi$. 

The existence of a compact global attractor for \eqref{eq-attrac} is well known
for the Sobolev subcritical case $p<3$. The first proofs in this case go back to
1985 
(\cite{Hale-book} and \cite{Haraux}), see \cite{Raugel} for other references.
The case $p=3$ as been studied in \cite{BV} and \cite{ACH}. For $p\in (3,5)$,
Kapitanski proved in \cite{Kapitanski} the existence of a compact global
attractor for \eqref{eq-attrac} if $\Omega$ is a compact manifold without
boundary and if $\gamma(x)=\gamma$ is a constant damping. Using the same
arguments as in the proof of our main result, we can partially deal with the
case $p\in (3,5)$ with a localized damping $\gamma(x)$ and with unbounded
manifold with boundaries. 

Assume that $f$ satisfies the assumptions of Theorem \ref{th-attrac}. 
We first claim that we can assume in addition that $f(x,0)=0$ on
$\partial\Omega$ in order to guarantee the Dirichlet boundary condition for
$f(x,u)$ if $u\in H^1_0(\Omega)$. Indeed, let $\varphi$ be the solution of $\Delta\varphi-\beta \varphi=f(x,0)$ with Dirichlet boundary conditions, which is well defined and smooth, since by assumptions $f(x,0)$ is smooth and compactly supported. Let
$\chi$ be a smooth compactly supported cut-off function such that $\Drond \chi\nu =0$ on $\partial
\Omega$ and $\chi\equiv 1$ in the ball $B(x_0,R)$ of Hypothesis \eqref{hyp-f3}. If $f(x,0)\neq 0$ on
$\partial\Omega$, we consider the new variable $\tilde u=u-\chi\varphi$ and the new equation
\begin{equation}\label{eq-attrac-new}
\partial^2_{tt}\tilde u + \gamma(x) \partial_t \tilde u = \Delta \tilde u - \beta \tilde u - \tilde f(x,\tilde u)
\end{equation}
where $\tilde f(x,\tilde u)=f(x,\tilde u+\chi\varphi(x))-\chi f(x,0)+2\grad\varphi.\grad\chi+\varphi\Delta\chi$. One directly checks that $u$ solves \eqref{eq-attrac} if and only if $\tilde u$ solves \eqref{eq-attrac-new}. Moreover $\tilde f$ is smooth, analytic in $\tilde u$ and satisfies Hypotheses \eqref{hyp-f2} and \eqref{hyp-f3} and in addition $\tilde f(x,0)=0$ on the boundary. Clearly, since our change of variables is a simple translation, 
obtaining a compact global attractor for \eqref{eq-attrac-new} is equivalent to
obtaining a compact global attractor for \eqref{eq-attrac}. In what follows, we
may thus assume that $f(x,0)=0$ in addition to the hypotheses of Theorem
\ref{th-attrac}.

The arguments of this paper show the following properties.
\begin{enum2}
\item {\bf The positive trajectories of bounded sets are bounded}. Indeed, 
\eqref{hyp-f3} implies that, for $x\not\in
B(x_0,R)$, we have $V(x,u)=\int_0^u f(x,\xi)d\xi\geq 0$. Moreover, for $x\in
B(x_0,R)$, $V(x,.)$ is non-increasing on $(-\infty,-R)$ and non-decreasing on
$(R,\infty)$. Thus, $V(x,u)$ is bounded from below for $x\in
B(x_0,R)$ and 
$$\forall (u,v)\in X~,~~E(u,v)\geq \frac 12 \|(u,v)\|_X^2 + vol(B(x_0,R)) \inf  
V~.$$
The Sobolev embeddings $H^1(\Omega)\hookrightarrow L^{p+1}(\Omega)$
shows that the bounded sets of $X$ have a bounded energy. Since the energy $E$
is non increasing along the trajectories of \eqref{eq-attrac}, we get that the
trajectory of a bounded set is bounded.

\item {\bf The dynamical system is asymptotically smooth}. The asymptotic
compactness exactly corresponds to the statement of Proposition
\ref{prop-asympt-reg}. Let us briefly explain why it can be extended to the
case where $f$ depends on $x$. The key point is the extension of Corollary
\ref{coro-DLZ}. First notice that, since we have assumed that $f(x,0)=0$ on $\partial\Omega$, $f(x,u)$ satisfies Dirichlet boundary condition if $u$ does. Then, it is not difficult to see that the discussion following
Theorem \ref{th-DLZ} can be extended to the case $f$ depending on $x$ by using 
estimates \eqref{hyp-f2}. Corollary \ref{coro-DLZ} follows then, except for a
small change: since it is possible that $f(x,0)\neq 0$ for some $x\in\Omega$,
the conclusion of Corollary \ref{coro-DLZ} should be replaced by 
$$\|f(x,v)\|_{L^1([0,T],H^{s+\varepsilon}_0(\Omega))}\leq C(1+
\|v\|_{L^\infty([0,T],H^{s+1}(\Omega)\cap H^1_0(\Omega))})~.$$
Then the proof of Proposition \ref{prop-asympt-reg} is based on Corollary
\ref{coro-DLZ}, the boundedness of the positive trajectories of bounded sets
(both could be
extended to the case where $f$ depends on $x$ as noticed above) and an
application of Proposition \ref{prop-positive-damping} outside of a large ball.
We conclude by noticing that, for $x$ large, $f(x,u)u\geq 0$ and
$\gamma(x)\geq
\alpha>0$ and thus Proposition \ref{prop-positive-damping} can still be applied
exactly as in the proof of Proposition \ref{prop-asympt-reg}.

\item {\bf The dynamical system generated by \eqref{eq-attrac} is gradient},
that is 
that the energy $E$ is non-increasing in time and is constant on a trajectory
$u$ if and only if $u$ is an equilibrium point of \eqref{eq-attrac}.
This last property is shown in Proposition \ref{prop-eq} for $f$ independent of
$x$ but can be easily generalised for $f=f(x,u)$. Notice that the proof of this
property is the one, where the analyticity of $f$ is required since the unique
continuation property of Section \ref{sect-UCP} is used. 
Finally, we enhance that the gradient structure of \eqref{eq-attrac} is
interesting from the dynamical point of view since it implies that any
trajectory $u(t)$ converges when $t$ goes to $+\infty$ to the set of equilibrium
points. 

\item {\bf The set of equilibrium points is bounded}.  The argument is similar
to the one of Corollary \ref{coro-uniNL}: if $e$ is an equilibrium point of
\eqref{eq-attrac} then \eqref{hyp-f3} implies that 
$$\int \frac12|\grad e|^2+\beta|e|^2 = -\int_\Omega f(x,e)e dx \leq -
vol(B(x_0,R)) \inf_{(x,u)\in\overline \Omega\times\RR} f(x,u)u~,$$
where we have bounded $f(x,u)u$ from below exactly as we have done for
$V(x,u)$ in i).
\end{enum2}
It is well known (see \cite{Hale-book} or Theorem 4.6 of \cite{Raugel} for examples) that
Properties i)-iv) yield the existence of a compact global attractor. Hence, we
obtain the conclusion of Theorem \ref{th-attrac}.

%%%%%%%%%%%%%%%%%%%%%%%%%%%%%%%%%%%%%%%%%%%%%%%%%%%%%%%%%%%%%%%%%%%%%%%%%%%%
%%%%%%%%%%%%%%%%%%%%%%%%%%%%%%%%%%%%%%%%%%%%%%%%%%%%%%%%%%%%%%%%%%%%%%%%%%%%

\end{document}